\documentclass[10pt]{article}
\usepackage{amsmath,enumerate,makeidx}
\usepackage{amssymb,geometry,latexsym,theorem}
\usepackage[all]{xy}
\usepackage[ps2pdf, bookmarks=true,
            bookmarksnumbered=true, breaklinks=true,
            pdfstartview=FitH, hyperfigures=false,
            plainpages=false, naturalnames=true,
            colorlinks=true,
            pdfpagelabels]{hyperref}
\usepackage[usenames]{color}
{\theorembodyfont{\rmfamily}\newtheorem{example}{Example}[section]}
{\theorembodyfont{\rmfamily}\newtheorem{definition}[example]{Definition}}
{\theorembodyfont{\rmfamily}\newtheorem{Definition}[example]{Definition}}
{\theorembodyfont{\rmfamily}}
\newtheorem{proposition}[example]{Proposition}
\newtheorem{theorem}[example]{Theorem}
{\theorembodyfont{\rmfamily}\newtheorem{remark}[example]{Remark}}
\newtheorem{corollary}[example]{Corollary}
\newtheorem{lemma}[example]{Lemma}

\newtheorem{blank}[example]{}
\newenvironment{proof}{\noindent {\bf Proof}}{\rule{0em}{1ex} %\phantom{z}
 \hfill $\Box$}

\def\Crs{\mathsf{Crs}}
\def\CRS{\mathsf{CRS}}
\def\Chn{\mathsf{Chn}}
\def\FTop{\mathsf{FTop}}
\def\FTOP{\mathsf{FTOP}}
\def\Top{\mathsf{Top}}

\def\Simp{\mathsf{Simp}}

\def\Sets{\mathsf{Sets}}

\def\Mod{\mathsf{Mod}}

\newcommand{\labra}[1]{\stackrel{#1}{\longrightarrow}}
\def\im{\mathrm{Im \,}}
\def\Cok{\operatorname{Coker}}
\def\Ker{\operatorname{Ker}}
\def\Ob{\operatorname{Ob}}
\def\px{{\rule{0em}{1.7ex}}}

\def\I{{\mathcal I}}

\def\C{\mathsf{C}}

\def\ge{\geqslant}
\def\leq{\leqslant}
\def\geq{\geqslant}

\def\subs{\subseteq}
\def\Cone{\operatorname{Cone}}
\def\bDelta{\boldsymbol{\Delta}}
\def\bUpsilon{\boldsymbol{\Upsilon}}
\def\Cyl{\operatorname{Cyl}}
\def\aI{I^\to}
\def\aZ{{\mathbb Z}^\to}

\def\eps{\varepsilon}
\def\subset{\subseteq}
\def\lan{\langle}
\def\ran{\rangle}
\def\Cone{\operatorname{Cone}}
\def\epb{\overline{\eps}}
\def\bS{\mathbb{S}}
\def\F{\mathbb{F}}

\def\colim{\operatorname{colim}}
\def\Piu{{\Pi^\Upsilon}}
\def\xybiglabels{\def\labelstyle{\textstyle}}

\begin{document}

\title{Normalisation for the fundamental crossed complex \\ of a simplicial set }
\author{Ronald Brown\thanks{Brown was supported by a Leverhulme Emeritus Fellowship, 2002-2004,
which also supported collaboration of the authors.}\\
r.brown@bangor.ac.uk\\
University of Wales, Bangor,
Dean St., \\
Bangor\\
Gwynedd LL57 1UT\\
U.K.\and Rafael Sivera\\{Rafael.Sivera@uv.es}\\
Departamento de {Geometr\'{i}a y Topolog\'{i}a}\\
          Universitat de Val\`encia,\\ 46100 Burjassot, Valencia\\
         Spain}

\maketitle
%\classification{18D10, 18G30, 18G50, 20L05, 55N10, 55N25, 55U10,
%55U99 }

%\keywords{crossed complex, simplicial set, normalisation, homotopy
%addition lemma, groupoid, crossed module, fundamental crossed
%complex}

\begin{center}
  {\it Dedicated to the memory of Saunders Mac Lane}
\end{center}
\begin{abstract}Crossed complexes are shown to have an algebra
sufficiently rich to model the geometric inductive definition of
simplices, and so to give a purely algebraic proof of the  Homotopy
Addition Lemma (HAL) for the boundary of a simplex. This leads to
the {\it fundamental crossed complex} of a simplicial set.  The main
result is a normalisation theorem for this fundamental crossed
complex, analogous to the usual theorem for simplicial abelian
groups, but more complicated to set up and prove, because of the
complications of the HAL and of the notion of homotopies for crossed
complexes. We start with some historical background, {and give a
survey of the required basic facts on crossed complexes.}
\footnote{MATHSCICLASS: 18D10, 18G30, 18G50, 20L05, 55N10, 55N25,
55U10,55U99
 }
\end{abstract}

\tableofcontents
\section*{Introduction} \label{sec:intro}
\addcontentsline{toc}{section}{Introduction} {\it Crossed complexes}
are analogues of chain complexes but with nonabelian features in
dimensions 1 and 2.  So one aim of the use of crossed complexes is
to increase the power of methods analogous to those of chain
complexes.

Crossed complexes  can incorporate information on presentations of
groups, or groupoids. Thus another aim is to bring features of the
fundamental group nearer to the centre of the toolkit of algebraic
topology. From the early 20th century,  the fundamental group
$\pi_1(X,p)$ has played an anomalous r\^ole in algebraic topology.
This invariant of a pointed space has the properties:
\begin{itemize}
  \item  it is nonabelian;
  \item it can under some circumstances be calculated precisely  by a
van Kampen theorem;
\item presentations of it are important;
\item it models all  pointed, connected,  weak homotopy 1-types.
\end{itemize}
Yet higher dimensional tools (homology groups, homotopy groups) were
generally abelian. An exception was the use of crossed modules, as
developed by J.H.C. Whitehead in \cite{jhcw:CHII}, and independently
by Peiffer, \cite{Peiff}, and Reidemeister, \cite{Reid}. They were
shown in \cite{maclane&jhcw} to model all pointed, connected, weak
homotopy 2-types (there called 3-types).

A possible resolution of this anomaly -- nonabelian invariant in
dimension 1, abelian invariants in dimensions greater than 1 --  has
appeared with the transition from groups to groupoids. The
fundamental groupoid $\pi_1(X,A)$ on a set $A$ of base points was
shown to have computational and conceptual utility,
\cite{B67,brownbook:2}. Groupoids were found to have 2-dimensional
nonabelian generalisations, for example crossed modules and forms of
double groupoids,  again with computational and conceptual utility,
\cite{BS1,BS76}. Whitehead's crossed modules derived from homotopy
theory could in some cases be {\it  calculated precisely as
nonabelian structures} by a 2-dimensional van Kampen type theorem,
whose proof used a relative homotopy  double groupoid,
\cite{bh1978,br:w}.

These results were generalised in \cite{bh:colimits} to all
dimensions,  using crossed complexes, whose definition and basic
theorems we recall in section \ref{sec:basic}. It is important that
these results are proved by working directly with homotopically
defined functors and without the use of traditional tools such as
homology and simplicial approximation.  Surveys of the use of
crossed complexes are in \cite{brown:hha,brown:fields}. We also
mention the work of Huebschmann, \cite{Hueb}, on the representation
of group cohomology by crossed $n$-fold extensions, and the use in
\cite{sharko} of Whitehead's methods from CHII to model
algebraically  filtrations of a manifold given by Morse functions.

We can also see  from the use of the Homotopy Addition Lemma in
proofs of the Hurewicz and relative Hurewicz Theorems,
\cite{GWWhitehead},  that the category $\Crs$ lies at the transition
between homology and homotopy. For further work on crossed
complexes, see for example
\cite{baues1,Baues-Cond,bauesT1,emmathesis,tonksthesis}. The works
with Baues refer to crossed complexes as `crossed chain complexes',
and, as with Huebschmann,  are in the one vertex case. It can be
argued that the category $\Crs$ gives a linear approximation to
homotopy theory: that is, crossed complexes can incorporate
presentations of the fundamental group(oid)  and its actions. They
do not incorporate, say, higher dimensional Whitehead products, or
composition operators. The tensor product of crossed complexes (see
later in \ref{thm:monclosed}) allows, analogously to work on chain
complexes, for  corresponding notions of an `algebra' and so for the
modelling of more structure, as in
\cite{brown:gilbert,Baues-Cond,bauesT1}.

The main result of this paper, Theorem  \ref{thm:normalisation},
extends the theory  by proving for crossed complexes an analogue of
a basic normalisation theorem for the traditional chain complex
associated with a simplicial set, due originally to Eilenberg and
Mac Lane, \cite{eil-mac}. We use the Homotopy Addition Lemma to
construct two crossed complexes associated to a simplicial set $K$,
the unnormalised $\Piu K$ and the normalised $\Pi K$, of which the
first is free on all simplices of $K$ and the second is free on the
nondegenerate simplices. It is important to have both crossed
complexes for applications of acyclic model theory; this is
analogous to the application of acyclic models to the usual singular
chain complex  by Eilenberg and Mac Lane in
\cite{eil-mac-homhot,eil-mac-acyclic}.  Indeed it was this relation
with acyclic model theory which motivated this investigation, and
which we plan to deal with elsewhere. As an example of a corollary
from our normalisation theorem  and properties of crossed complexes
we obtain the well known fact that the projection $\|K \| \to |K|$,
from the thick to the standard geometric realisation of a simplicial
set $K$, is a homotopy equivalence (Corollary
\ref{cor:thicktothin}).

{To make this paper self-contained, we give a fairly full account of
the necessary properties of crossed complexes, so that this paper
can form an introduction to their use. The structure of this paper
is as follows: Section \ref{sec:hal} gives an introduction to the
Homotopy Addition Lemma: this Lemma is essential for describing the
fundamental crossed complex of a simplicial set. Section
\ref{sec:basic} gives an account of the basic results on crossed
complexes that are needed. Section \ref{sec:chains} states the main
relation with chain complexes with a groupoid of operators; this is
useful for understanding many constructions on crossed complexes.
Sections \ref{sec:generating}, \ref{sec:normal}, \ref{sec:free} give
brief accounts of generating crossed complexes, normal subcrossed
complexes and free crossed complexes;  the first two topics are not
easily available in the literature. } Sections
\ref{sec:cyl-homotopy}, \ref{sec:cones-HAL} give  specific rules for
the crossed complex constructions of cylinders and homotopies,  and
then cones, and so allow the algebraic deduction of  the Homotopy
Addition Lemma. Section \ref{sec:PiK} defines the (non normalised)
fundamental crossed complex of a simplicial set. Section
\ref{sec:0-normalisation} normalises this crossed complex at
$\eps_0$. Section \ref{sec:normalisation} gives the full
normalisation theorem.

\noindent {\bf Acknowledgement} We would like to  thank Philip
Higgins for a long and happy collaboration on the background to this
paper, and a referee for helpful comments.

\section{The Homotopy Addition Lemma} \label{sec:hal}
The normalisation theorem for simplicial abelian groups (see, for
example, Eilenberg-Mac~Lane \cite[Theorem 4.1]{eil-mac},  Mac~Lane
\cite[\S VIII.6]{maclane:hom}), is of importance in homological
algebra and in geometric applications of simplicial theory. It is
based on the formula, fundamental in much of simplicial based
algebraic topology and homological algebra, that if $x$ has
dimension $n$ then
\begin{equation} \label{eq:abelianboundary}
\partial x=
\sum_{i=0}^n (-1)^i
\partial_i x,\end{equation}
which can be interpreted intuitively as: `the boundary of a simplex
is the alternating sum of its faces'. The setting for this formula
is the theory of chain complexes: these are sequences of morphisms
of abelian groups (or $R$-modules) $\partial: A_n \to A_{n-1}$ such
that $\partial \partial=0$. {For simplicial abelian groups}, each
$\partial_i$ is a morphism of abelian groups and the formula
\eqref{eq:abelianboundary} is just the alternating sum of morphisms.
Thus we have a chain complex $(A, \partial)$. Further, if $(DA)_n$
is for $n \geq 0$ the subgroup of $A_n$ generated by degenerate
elements, then $DA$ is a contractible subchain complex of
$(A,\partial)$. This is the normalisation theorem.

In {\it homotopy}, rather than homology,  theory, there is another
and more complicated basic formula,  known as the {\it Homotopy
Addition Lemma} (HAL) (or theorem) \cite[\S 12]{Blakers},
\cite{Hu-hal}, \cite[Theorem IV-6.1, p174 ff.]{GWWhitehead}. The
intuition of `the boundary of a simplex' given in
\eqref{eq:abelianboundary} is strengthened in the HAL by taking
account also of:
\begin{itemize}
\item a {\it set} of base
points (the vertices of the simplex); \item nonabelian structures in
dimensions 1 and 2; \item operators of dimension 1 on dimensions
$\geq 2$.
\end{itemize}From our standpoint, the set of base points is
taken into account through the use of groupoids in dimension 1,
while the boundary from dimension 3 to dimension 2 uses crossed
modules of groupoids. This leads to basic formulae, which with our
conventions are as follows:

In dimension 2 we have a groupoid rule:
\begin{equation}\label{HAL2}
\delta_2x =  -\partial _1 x + \partial_2 x +
\partial_0 x, \tag{HAL2}
\end{equation}
which is represented by the diagram
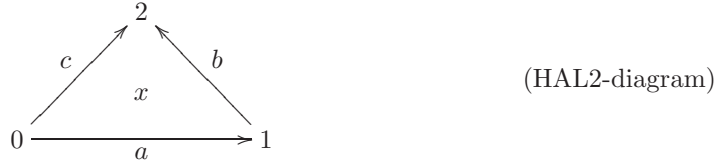
\begin{equation}\label{HAL2-diagram}
\vcenter{\xybiglabels\xymatrix@=3pc{& 2& \\
0 \ar [ur] ^c \ar [rr] _a & \ar @{}[u]|(0.35)x & 1 \ar [ul] _ b}}
\tag{HAL2-diagram}
\end{equation}
and the easy to understand formula \eqref{HAL2} says that $\delta_2x
=-c +a +b$. Note that we use additive notation throughout for group
or groupoid composition. Our  convention is that  the {\it base
point} of an ordered $n$-simplex $x$  is the final vertex
$\partial_0^n x$.

In dimension 3 we have the nonabelian rule:
\begin{equation}\label{eq:HAL3}
\delta_3 x = (\partial_3x) ^{\partial^2_0x} - \partial_0 x
-\partial_2 x +\partial_1x. \tag{HAL3}\end{equation} Understanding
of this is helped by considering the diagram
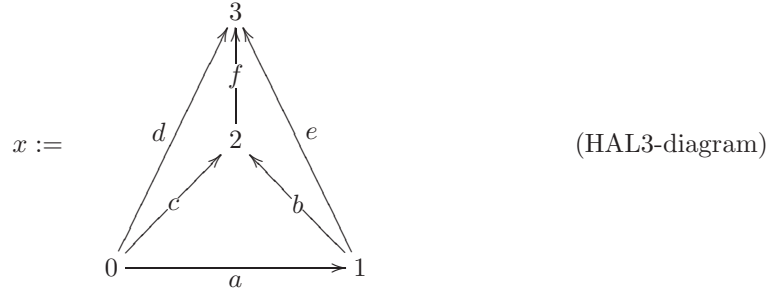
\begin{equation} \label{HAL3-diagram}x:=\quad
\xybiglabels\vcenter{\xymatrix@=3pc{&3& \\
& 2 \ar [u]|{\rule[-0.3ex]{0em}{2ex}f} & \\
0 \ar[uur]^d \ar [ur] |c \ar [rr] _a && 1 \ar [ul] |b \ar [uul] _e
}}\tag{HAL3-diagram}\end{equation} in which  $f=\partial_0^2x$. The
base point of the above 3-simplex $x$ is $3 $, while the base point
of $\partial_3 x$ is 2. The exponent $f$ relocates $\partial_3 x$ to
have base point at 3, and so  yields a well defined formula. Given
the labelling in \eqref{HAL3-diagram} we have the groupoid formula
$$-f+(-c+a +b)+f -(-e+b+f) -(-d+a+e)+(-d+c+f)=0. $$
This is a translation of the rule $\delta_2 \delta_3=0$, provided we
assume $\delta_2(y^f)=-f + \delta_2 y +f$, which is the first rule
for a crossed module.

In dimension $n\ge 4$ we have the abelian rule  with operators:
\begin{equation}\label{HALge-4}
 \delta_nx  =  (\partial_n
x)^{u_n x   } + \sum _{i=0}^{n-1} (-1)^{n-i}
\partial_ix \quad \mbox{ for $n \ge 4$},\tag{HAL$\ge 4$}
\end{equation}
where $u_nx= \partial_0^{n-1} x$. We have some difficulty in drawing
a diagram for this! These,  or analogous, formulae underly much
nonabelian cohomology theory.

The rule $\delta_{n-1}\delta _n=0$ is  straightforward to verify for
$n>4$, through working in abelian groups;  for $n=4$ we require the
second crossed module rule, that for $x,y$ of dimension 2
$$-y+x+y= x^{\delta_2 y}.$$ A consequence is that $\Ker \delta_2$,
is central. Hence so also is $\im \delta_3$, since we have verified
$\delta_2\delta_3=0$. The type of argument that is used for the
$n=4$ case, \cite{GWWhitehead}, and which we shall use later, is the
simple:
\begin{lemma}\label{lem:centrality}
  If $\gamma=\alpha +\beta $   is  a central element in a group, then $\gamma=\beta+\alpha$.
\end{lemma}\begin{proof}
\begin{equation}
\gamma - \alpha - \beta =-\alpha +\gamma -\beta=0. \tag{centrality}
\end{equation}
\end{proof}

The above formulae are not exactly as will be found generally in the
literature; they follow from our conventions given in section
\ref{sec:cyl-homotopy} for crossed complexes and their `cylinder
object' $\I \otimes C$ and cone object $\Cone(C)$. Thus  the chain
complex and homological boundary formula \eqref{eq:abelianboundary}
becomes a much more complicated  result  in homotopy theory, as the
formulae in the Homotopy Addition Lemma. Formulae of this type occur
frequently in mathematics, for example in the cohomology of groups,
\cite{BRazak:LMS99},  in differential geometry, \cite{kock}, and in
the cohomology of stacks, \cite{breen}.

The formal structure required for this HAL is known as a {\it
crossed complex of groupoids}, and such structures  form the objects
of a category which we write $\Crs$.

We give more details of this category in section \ref{sec:basic}. It
is complete and cocomplete. It contains `free' objects, satisfying
the universal property that a morphism $f: F \to C$ from a free
crossed complex is defined by its values on a free basis, subject to
certain geometric conditions. Note that  in the formulae for the
HAL, the $\partial_i$ are {\it not morphisms}, but $x$ and all
$\partial_i x$ are elements of a free basis.

Our first aim here is to show how the HAL for a simplex fits neatly
into an algebraic pattern in crossed complexes, using a cone
construction $\Cone(B)$ for a crossed complex $B$. We define
algebraically and inductively an `algebraic' or `crossed complex
simplex' $a\Delta^n$ by
\begin{equation}
a\Delta^0 =\{0\}, \quad a\Delta^{n} = \Cone(a\Delta^{n-1}).
\end{equation}
Our  conventions for the tensor product, \cite{bh:tens}, ensure that
this definition yields algebraically {precisely} the HAL given
above.

{Our main result is an application to the (non normalised)
fundamental crossed complex, written here  $\Pi^\Upsilon K$, of a
simplicial set $K$. This is defined to be the free crossed complex
on the elements of $K_n, n \ge 0$, with boundary given by the
homotopy addition lemma HAL. Thus $\Pi^\Upsilon K$ contains basis
elements which are degenerate simplices, of the form $\eps_i y$ for
some $y$. Full details are given in section \ref{sec:PiK}.

We may also construct $\Pi^\Upsilon K$ as a coend as follows. Let
$\bDelta$ be the simplicial operator category, so that a simplicial
set $K$ is a functor $\bDelta^{op} \to \Sets$. Let $\bUpsilon$ be
the subcategory of $\bDelta$ generated by the injective maps, i.e.
those which correspond to simplicial face operators. Then we can see
the unnormalised fundamental crossed complex of $K$ as the coend in
the category of crossed complexes
$$\Pi^\Upsilon K= \int^{\bUpsilon,n} K_n \times a\Delta^n.  $$
We are also interested in the {\it normalised } crossed complex,
defined as the coend
$$\Pi K= \int^{\bDelta,n} K_n \times a\Delta^n.  $$
In Section \ref{sec:normalisation} we complete the proof of:
\begin{theorem}[Normalisation theorem]
For a simplicial set $K$, the quotient morphism $p: \Piu K \to \Pi
K$  is a homotopy equivalence with section $q: \Pi  K \to \Piu K$,
and the quotient crossed complex $\Pi  K$ is free.
\end{theorem}}
This has application to the usual thick and standard geometric
realisations of a simplicial set $K$, defined respectively as
coends: \begin{align*}
 \| K\| &= \int^{\bUpsilon,n} K_n \times \Delta^n, \\
|K| &= \int^{\bDelta,n} K_n \times \Delta^n,
\end{align*}
where $\Delta^n$ is the geometric simplex. Then the normalisation
theorem together with standard properties of crossed complexes,
implies:
\begin{corollary}\label{cor:thicktothin}
For a simplicial set $K$, the projection $\|K \| \to |K|$ from the
thick to the standard geometric realisation is a homotopy
equivalence.
\end{corollary}
\begin{proof}
We use the Higher Homotopy  van Kampen Theorem (HHvKT) of
\cite{bh:colimits} to give natural isomorphisms $$\Piu K \cong \Pi
(\|K \|_*), \Pi K \cong \Pi (|K|_*).$$  It is immediate that the
projection induces an isomorphism of fundamental groupoids and of
the homologies of the universal covers at all base points.
\end{proof}

We assume work on groupoids as in \cite{brownbook:2,Higginsbook}.

\section{Basics on  crossed complexes} \label{sec:basic}
Crossed complexes, first called  {\it group systems},  were first
defined, in the one vertex case, in 1946 by Blakers in
\cite{Blakers}, following a suggestion of Eilenberg. He combined
into a single structure the fundamental group {$\pi_1(X,p)$} and the
relative homotopy groups {$\pi_n(X_n, X_{n-1},p), n \geq 2$}
associated to a reduced filtered space $X_*$, i.e. when $X_0$ is a
singleton {$\{p\}$}. We now call this structure the {\it fundamental
crossed complex} $\Pi(X_*)$ of the filtered space (see below).

Blakers' concept was taken up in J.H.C. Whitehead's deep paper
`Combinatorial homotopy theory II' (CHII) \cite{jhcw:CHII}, in the
reduced and free case, under the term `homotopy system'; this paper
is much less read than the previous paper `Combinatorial homotopy I'
(CHI) \cite{jhcw:CHI}, which introduced the basic concept of
$CW$-complex. We give below a full definition of the category $\Crs$
of crossed complexes: our viewpoint, following that of CHII,  is
that $\Crs$ should be seen as a basic  category for applications in
algebraic topology, with better realisability properties,
\cite{bh:class91,brmopw}, than the more usual chain complexes with a
group of operators, \cite{bh:chn}.

We  use relative homotopy theory to construct the functor
\begin{equation} \Pi: \FTop \to \Crs
\end{equation}
where $\FTop$ is the category of {\it filtered spaces}, whose
objects
$$X_*: X_0 \subs X_1 \subs \cdots \subs X_n \subs \cdots \subs
X_\infty $$ consist of a compactly generated topological space
$X_\infty$ and an increasing sequence of subspaces $X_n, n \ge 0$.
The morphisms $f: X_* \to Y_*$ of $\FTop$ are maps $f: X_\infty \to
Y_ \infty$ such that for all $n \ge 0$ $f(X_n ) \subs Y_n$.

The functor $\Pi$ is given on a filtered space $X_*$ by
\begin{equation}
  (\Pi X_*) _n = \begin{cases}
   X_0 & \text{ if } n=0,\\
   \pi_1(X_1,X_0) & \text{ if } n=1,\\
   \pi_n(X_n,X_{n-1},X_0) & \text{ if } n\ge 2.
  \end{cases}
\end{equation}
Here $\pi_1(X_1,X_0)$ is the fundamental groupoid of $X_1$ on the
set $X_0$ of base points, and $\pi_n(X_n,X_{n-1},X_0)$ is the family
of relative homotopy groups {$\pi_n(X_n, X_{n-1},p)$} for all {$p
\in X_0$}.

If we write $C_n= (\Pi X_*)_n$, then we find that there is a
structure of a family of groupoids over $C_0$ with source and target
maps $s,t$: \begin{equation}\label{eq:crossedcompl} \xymatrix{
\cdots \ar [r] & C_n \ar @<-.05ex>[d] ^{t}\ar@<0.2ex> [r]
^-{\delta_n} & C_{n-1} \ar@<0.2ex> [r]   \ar @<-1.2ex>[d]  ^{t} &
\cdots \ar@<0.2ex>  [r] &
 C_2\ar @<0.2ex> [r]  ^-{\delta_2} \ar @<-.05ex>[d]  ^{t} & C_1
\ar @<0.0ex>[d]  ^(0.45){t} \ar @<-1ex>[d]  _(0.45){s} \\ &
C_0&C_0\rule{0.5em}{0ex}  & & \rule{0.5em}{0ex} C_0 &
\rule{0em}{0ex} C_0}\end{equation}  in which: for $n \geq 2$, $C_n$
is totally disconnected, i.e. $s=t$;  $C_1$ operates (on the right)
on $C_n, n \ge 2$, and on the family of vertex groups of $C_1$ by
conjugation; and the axioms are, in addition to the usual operation
rules, that:\begin{enumerate}[CC1)]
\item  $s\delta_2 =t\delta_2$, $\delta_{n-1}\delta_n = 0$;\\
\item $  \delta_n$ is an operator morphism;\\
\item $\delta_2: C_2 \to C_1$ is a crossed module; \\
\item for $n \ge 3$, $C_n$ is abelian  and $\delta_2 C_2$ operates
trivially on $C_n$.\end{enumerate} It will be convenient to write
all group and groupoid compositions additively, and the operations
as $x^a$: if $a:p \to q$ in dimension 1, then $p=sa, q=ta$, and if
further $b: q \to r$ then $a+b: p \to r$; if  $ n \geq 2$ and $x \in
C_n(p)$, then $tx=p$ and $x^a \in C_n(q)$.

These laws  CC1)-CC4) for $C= \Pi X_*$ reflect basic facts in
relative homotopy theory. Indeed, for any crossed complex $C$ there
is a filtered space $X_*$ such that $C \cong \Pi X_*$,
\cite{bh:colimits}. These laws also define the objects of the
category $\Crs$ of crossed complexes. The morphisms $f: C \to D$ of
crossed complexes consist of groupoid morphisms $f: C_n \to D_n, n
\ge 1$, preserving all the structure.

A crossed complex $C$ has a {\it fundamental groupoid} $\pi_1 C$
defined to be $C_1/(\delta_2 C_2)$, whose set of components is
written $\pi_0 C$, and also called the {\it set of components} of
$C$. It also has a family of {\it homology groups} given for $n \geq
2$ by
\begin{equation*} H_n(C,p)= \Ker(\delta_n: C_n(p) \to C_{n-1}(p)) /
\im(\delta_{n+1}( C_{n+1}(p) \to C_n(p)),
\end{equation*}which can  be seen to be a module over $\pi_1 C$.
A morphism $f:C \to D$ of crossed complexes induces a morphism of
fundamental groupoids and homology groupoids, and is called a {\it
weak equivalence} if it induces an equivalence of fundamental
groupoids and isomorphisms $H_*(C,p) \to H_*(D,fp)$ for all $p \in
C_0$.

If $X_*$ is the skeletal filtration of a $CW$-complex $X$, then
{(see \cite{jhcw:CHII}) } there are natural isomorphisms
\begin{equation*} \pi_1(\Pi X_*) \cong \pi_1(X, X_0), \qquad H_n(\Pi
X_*, p) \cong H_n (\widetilde{X}_p),
\end{equation*}
where $\widetilde{X}_p$ is the universal cover of $X$ based at $p$.
It follows from  this and Whitehead's theorem {from \cite{jhcw:CHI}}
that if $f: X \to Y$ is a cellular map of $CW$-complexes $X,Y$ which
induces a weak equivalence $\Pi f: \Pi X_* \to \Pi Y_*$, then $f$ is
a homotopy equivalence.

The following additional facts on crossed complexes were found in a
sequence of papers by Brown and Higgins:
\begin{blank}[\cite{bh:colimits} ]\label{thm:colim} The functor
 $\Pi: \FTop \to \Crs$
from the category of filtered spaces to crossed complexes preserves
certain colimits.
\end{blank}
\begin{blank}[\cite{bh:tens}]\label{thm:monclosed}
The category $\Crs$ is monoidal closed, with an exponential law of
the form
\begin{equation} \label{equ:exponent}
  \Crs(A \otimes B, C) \cong \Crs(A, \CRS(B,C)). \tag{exponential law}
\end{equation}
\end{blank}
\begin{blank}[\cite{brown-gol}]\label{thm:hom}
The category $\Crs$ has a {\it unit interval object} written  $\{0\}
\rightrightarrows \I$, which is essentially just the indiscrete
groupoid on two objects $0,1$, and so has in dimension $1$ only one
element $\iota :0 \to 1$. For a crossed complex $B$, this gives rise
to a {\it cylinder object} \begin{align*} \Cyl(B) &= ( B \rightrightarrows \I \otimes B),%\\
%\intertext{and so a `cone construction'}
%  \Cone(B) &= \Cyl(B) /(\{1\} \otimes B).
\end{align*}
and so a homotopy theory for crossed complexes.
\end{blank}
The  result \ref{thm:colim} is a kind of Higher Homotopy van Kampen
Theorem (HHvKT). Among its consequences are the relative Hurewicz
Theorem, seen from this viewpoint as a relation between $\pi_n(X,A)$
and $\pi_n(X \cup CA)$. It also implies nonabelian results in
dimension 2 not obtained by other means, \cite{B:coproducts,br:w}.
{The proof of the HHvKT, \cite{bh:colimits}, uses cubical higher
homotopy groupoids and is independent of standard methods in
algebraic topology, such as homology or simplicial approximation.

The monoidal closed structure for crossed complexes allows us to
define homotopies for  morphisms $B \to C$ of crossed complexes as
morphisms $\I \otimes B \to C$, or, equivalently, as morphisms $\I
\to \CRS(B,C)$. The detailed structure of this  {\it cylinder object
$\I \otimes C$}, \cite{Kamps-Porter}, will be given in section
\ref{sec:cyl-homotopy}. The model structure of this homotopy theory
is developed in \cite{brown-gol}.

The full structure of the internal hom $\CRS(B,C)$  is quite
complicated. This complexity is also reflected in the structure of
the tensor product $A \otimes B$ of crossed complexes $A,B$: it is
generated in dimension $n$ by elements $a \otimes b$ where {$a \in
A_l, b \in B_k, l+k=n$}; the full list of structure and laws is
again quite complex (see \cite{bh:tens}).

We need that the  category $\FTop$  of filtered spaces is monoidal
closed with an exponential law
\begin{equation}
  \FTop(X_* \otimes Y_*, Z_*) \cong \FTop(X_*, \FTOP(Y_*,Z_*)).
\end{equation}
Here $(X_* \otimes Y_*)_n= \bigcup_{p+q=n} X_p \times Y_q$. A
standard example of a filtered space is a $CW$-complex with its
skeletal filtration, and among the $CW$-complexes we have the
$n$-ball $E^n$ with its cell structure
\begin{equation}
  E^n= \begin{cases}
    e^0 & \text{ if } n=0,\\
    e^0_\pm \cup e^1&  \text{ if } n=1, \\
    e^0 \cup e^{n-1} \cup e^n & \text{ if } n>1.
  \end{cases}
\end{equation}
The complications of the cell structure of $E^m \times E^n$ are
modelled in the tensor product of crossed complexes, as shown by:
\begin{blank}[\cite{bh:class91}]\label{thm:nattrans}
There is for filtered spaces $X_*,Y_*$ a natural transformation $$
\eta: \Pi X_* \otimes \Pi Y_* \to \Pi (X_* \otimes Y_*),$$  which is
an isomorphism if $X_*,Y_*$ are the skeletal filtrations of
$CW$-complexes, {\em \cite{bh:class91}}, and more widely,
\em{\cite{baues-br}}.
\end{blank}

\begin{remark}  From the early days, basic results  of relative homotopy
theory have been proved by relating the geometries of cells and
cubes. This geometric relation was translated  into a relation
between algebraic theories in several papers, particularly
\cite{bh:algcub,bh:tens}, which  give an equivalence of monoidal
closed categories between a category of `cubical $\omega$-groupoids'
and the category $\Crs$. While many constructions and proofs are
clearer in the former category,  both categories are required for
some results. For example the natural transformation $\eta$ of
\ref{thm:nattrans} is easy to see in the cubical category,
\cite{bh:class91}. For a survey on crossed complexes and their uses,
see \cite{brown:fields}. A book in preparation, \cite{Brown-Sivera},
is planned  to give a full account of all these main properties, and
make the theory more accessible and hopefully more usable.
\end{remark}

An important result on crossed complexes is the following:
\begin{blank}[\cite{bh:class91}]\label{thm:class}
There is a {\it classifying space functor} $B: \Crs \to \Top$ and a
homotopy classification theorem
$$[X, BC] \cong [\Pi X_*, C]$$for a $CW$-complex $X$ with its
skeletal filtration, and crossed complex $C$.
\end{blank}
This result is a homotopy classification theorem in the non simply
connected case, and includes many classical results. It relates to
Whitehead's comment in \cite{jhcw:CHII} that  crossed complexes have
better realisation properties than chain complexes with a group of
operators. It is also relevant to nonabelian cohomology, and
cohomology with local coefficients, as discussed in
\cite{bh:class91}. See also \cite{bullejosetal,
martins,martins-porter}.

\section{Crossed complexes and chain complexes}\label{sec:chains}
As is clear from the definition, crossed complexes differ from chain
complexes of modules  over groupoids only in dimensions 1 and 2. It
is useful to make this relationship more precise. We therefore
define a category $\Chn$ of such chain complexes which is monoidal
closed, give a functor $\nabla: \Crs \to \Chn$, and state that this
functor is monoidal and has a right adjoint. The results of this
section are taken from \cite{bh:chn}, which develops results from
\cite{jhcw:CHII}.  See also \cite{brohu} for the low dimensional and
reduced case.

We first define the category $\Mod$ of {\it modules over groupoids}.
We will often  write $G_0$ for $\Ob G$ for a groupoid $G$. The
objects of the category $\Mod$ are pairs $(G,M)$ where $G$ is a
groupoid and $M$ is a family $M(p), \, p \in G_0$, of disjoint
abelian groups on which $G$ operates. The notation for this will be
the same as for the operations of $C_1$ on $C_n$ for $n \geq 3$ in
the definition of a crossed complex. The morphisms of $\Mod$ are
pairs $(\theta,\phi): (G,M) \to (H,N)$ where $\theta: G \to H$ is a
morphism of groupoids and $\phi: M \to N$ is a family $\phi(p): M(p)
\to N(\theta p)$ of morphisms of abelian groups preserving the
operations. Instead of writing $(G,M)$ we often say $M$ is a
$G$-module. For a morphism $\phi: M \to N$ of $G$-modules it is
assumed that $\theta=1_G$.

A {\it chain complex $C$  over a groupoid} $G$ is a sequence of
morphisms of $G$-modules $\partial:C_n \to C_{n-1}, \; n \geq 1$,
such that $\partial\partial =0$. A {\it morphism of chain complexes}
consists of a morphism $\theta: G \to H$ of groupoids and a family
$\phi_n: C_n \to D_n$ such that $(\theta, \phi_n) $ is a morphism of
modules and also $\partial \phi_n= \phi_{n-1} \partial, n \geq 0$.
This defines the category $\Chn$.

The category $\Mod$ is monoidal closed. We define here only the
tensor product: for modules $(G,M),\; (H,N)$ we set $$(G,M) \otimes
(H,N) = (G \times H, M \otimes N)$$ using the product of groupoids
and with $$(M \otimes N)(p,q)= M(p) \otimes N(q),$$the usual tensor
product of abelian groups, and action the product action $(m\otimes
n)^{(g,h)}= m^g \otimes n^h$.

The category $\Chn$ is also monoidal closed with the usual tensor
product $(G,C) \otimes (H,D)= (G \times H, C \otimes D)$ where $(C
\otimes D)_n= \bigotimes_{p+q=n}C_p \otimes D_q$. For  details of
the internal hom, see \cite{bh:chn}.

A particular module we need is $(G, \aZ G)$, also written $\aZ G$,
for a groupoid $G$. If $p \in G_0$, then $\aZ G(p)$ is the free
abelian group on the elements of $G$ with final point $p$, and with
operations induced by the right action of $G$. Note that in contrast
to the single object case, i.e. of groups, we obtain a module and
not an analogue of a ring. A  set $J$ defines a discrete groupoid on
$J$ also written $J$ and so a module $(J, \aZ J)$: when the set $J$
is understood, we abbreviate this module to $\aZ$. The {\it
augmentation map} in this context is given as usual by the sum of
the coefficients. It is a   module morphism $\eps: (G, \aZ G) \to
(G_0, \aZ)$ and its kernel is the {\it augmentation module} $(G, \aI
G)$ which we abbreviate to $\aI G$.

Let $\psi: G \to H$ be a groupoid morphism which is bijective on
objects, and let $(H,N)$ be a module. We need the generalisation to
groupoids of the universal derivations of Crowell \cite{crowell}. A
$\psi$-derivation $d:G \to N$ assigns to each $g \in G$ with final
point $p$ an element $d(g) \in N(\psi p)$ satisfying the rule that
$d(g' + g)= d(g')^{\psi g}+ d(g)$ whenever $g'+g$ is defined in $G$.
The $\psi$-derivation $d$ is {\it universal} if given any other
$\psi$-derivation $d': G \to L$ where $L$ is an $H$-module, there is
a unique $H$-module morphism $f:L \to N$ such that $fd'=d$. The
construction of the universal $\psi$-derivation is straightforward,
and is written $\alpha: G \to D_\psi$.

\begin{blank} \label{thmII1:crossedtochain} Let $C$ be a crossed complex, and let
$\phi: C_1 \to G$ be a cokernel of $\delta_2$ of $C$. Then there are
$G$-morphisms
\begin{equation} \label{eq:exactchain}
C^{\rm ab}_2 \labra{\partial_2} D_{\phi} \labra{\partial_1 } \aZ G
\end{equation}
such that the diagram \begin{equation} \label{eq:crossedtochain}
 \xybiglabels \vcenter{\xymatrix@C=1.5pc{
\cdots \ar[r] & C_n \ar[d]_{=} \ar[r]^{\delta_n}& C_{n-1} \ar[d]_{=}
\ar[r] & \ldots \ar[r] & C_3 \ar[d]_{=} \ar[r]^{\delta_3}& C_2
\ar[d]^{\alpha_2}
\ar[r]^{\delta_2} & C_1 \ar[d]^{\alpha_1} \ar[r]^{\phi} & G \ar[d]^{\alpha_0} \\
\cdots \ar[r] & C_n \ar[r]_{\partial_n} & C_{n-1} \ar[r]& \cdots
\ar[r] & C_3  \ar[r]_{\partial_3} & C^{\rm ab}_2 \ar[r]_{\partial_2}
& D_{\phi} \rule[1ex]{0em}{1ex} \ar[r]_{\partial_1 } & \aZ G ^\px
\rule[-1ex] {0em}{1ex}}}
\end{equation} commutes and the lower line is a chain complex over
$G$. Here $\alpha_2$ is abelianisation, $\alpha_1$ is the universal
$\phi$-derivation, $\alpha_0$ is the $G$-derivation $x \mapsto
x-1_q$ for $x \in G(p,q)$, as a composition $G \to \aI G \to \aZ G$,
and $\partial_n = \delta_n$ for $n \geqslant 4$. Further
\begin{enumerate}[\rm (i)]
\item the sequence \eqref{eq:exactchain} is exact at $D_\phi$
and the image of $\partial_1$ is the augmentation module $\aI G$;
\item if $C_1$ is a free groupoid on a generating graph $X_1$ and
$\psi $ is surjective, then  $D_\phi$ is the free $H$-module on the
basis at $p \in H_0$ of elements $x$ of $X_1$ such that $\psi tx
=p$;
\item if $C_1$ is free, then $\alpha_2$ is injective on $\Ker
\delta_2$.
\end{enumerate}
\end{blank}

\begin{blank}\begin{enumerate}[\rm (i)]
\item The bottom row of diagram \eqref{eq:crossedtochain} defines a
functor $\nabla: \Crs \to \Chn$, which has a right adjoint. Hence
$\nabla$ preserves colimits.
\item The functor $\nabla$ preserves tensor products: there is a natural equivalence for crossed
complexes $A,B$ $$\nabla(A) \otimes \nabla (B) \cong \nabla(A
\otimes B).
$$
\end{enumerate}
\end{blank}
This last result shows that the major unusual complications of the
tensor product of crossed complexes occur in dimensions 1 and 2.
These cases are analysed in \cite{bh:tens}.
\begin{remark}
These results show the close relation of crossed complexes and these
chain complexes. The functor $\nabla$ loses information. Whitehead
remarks in [CHII] that (using our terminology)  these chain
complexes have less good realisation properties even than free
crossed complexes. Indeed, the problem of which 2-dimensional free
chain complexes are realisable by a crossed complex is known to be
hard.
\end{remark}

\section{Generating  {structures}}\label{sec:generating}
Let $C$ be a crossed complex, and let $R_*$ be a family of subsets
$R_n \subset C_n$ for all $n \geq 0$. We have to explain what is
meant by  the subcrossed complex $\langle R_* \rangle $ of $C$ {\it
generated} by $R_*$.

A formal definition of $B=\langle R_* \rangle $ is easy: it is the
smallest sub-crossed complex $D$ of $C$ such that $R_n \subset D_n$
for all $n \geq 0$, and so is also the intersection of all such $D$.
A direct construction is as follows. We set $$B_0= R_0 \cup sR_1
\cup \bigcup_{n \geq 1} tR_n.$$ Let $B_1$ be the subgroupoid of
$C_1$ generated by $R_1 \cup \delta_2(R_2)$ and the identities at
$B_0$. Let $B_2$ be the subcrossed $B_1$-module of $C_2$ generated
by $R_2$.  For $n \geq 3$, let $B_n$ be the sub-$B_1$-module of
$C_n$ generated by $R_n \cup \delta(R_{n+1}) $ and the identities at
elements of  $B_0$.

Note that this definition is inductive. The usual property of a
generating structure holds: thus if $R_*$ generates $C$, i.e. $\lan
R_* \ran = C$, and $f,g:C \to D$ are two crossed complex morphisms
which agree on $R_*$, then $f=g$. This is proved by induction.

We say the family $R_*$ is a {\em generating structure} for a
subcrossed complex $B$ of $C$ if for each $n >0$ the boundaries in
$C$ of elements of $R_n$ lie in the subcrossed complex generated by
the $R_i$ for $i < n$, and $R_*$ generates $B$.

\section{Normal subcrossed complexes}\label{sec:normal}
We assume work on normal subgroupoids, as in
\cite{brownbook:2,Higginsbook}.

\begin{definition}A subcrossed complex $A$ of a crossed complex $C$ is
called {\it normal} in $C$ if: \begin{enumerate}[N1)] \item $A$ is
wide in $C$, i.e. $A_0=C_0$; \item $A_1$ is a totally disconnected
normal subgroupoid of $C_1$; \item for $n \geq 2$, $A_n$ is
$C_1$-invariant, i.e. $a \in A_n$, {$x \in C_1$ and $a^x$ is defined
implies $a^x\in A_n$}; \item for $n \geq 2$, if  $a \in A_1$, {$x
\in C_n$, and $x^a$ is defined then $x-x^a \in A_n$}. \hfill $\Box$
\end{enumerate}
\end{definition}
Note that N3) implies that $A_2$ is a normal subgroupoid of $C_2$
since {$-x +a+x=a^{\delta_2 x}$}.

The above are necessary and sufficient conditions for $A$ to be the
kernel of a morphism $C \to D$ of crossed complexes {which is
injective on objects} for some $D$, as the following proposition
shows.
\begin{proposition}
  If $A$ is a normal subcrossed complex of the crossed complex $C$,
  then the family of quotients $C_n/A_n, \; n\ge 1$ inherits the structure of
  crossed complex, which we call the quotient crossed complex $C/A$.
\end{proposition}
{We leave the proof to the reader. }

Let $R_*$ be a family of subsets of the crossed complex $C$ as in
the previous section, and such that $R_1$ is totally disconnected,
i.e. just a family of subsets of vertex groups of the groupoid
$C_1$. We say that $R_*$ {\it normally generates} a subcrossed
complex $A$ of $C$ if $A$ is the smallest wide  normal subcrossed
complex of $C$ containing $R_*$, and then we say $A$ is the {\it
normal closure} of $R_*$ in $C$, and write $A=\lan \lan
R_*\ran\ran$. We also say $R_*$ is a {\it normal structure in $C$}
if for each $n>0$ the boundaries of elements of $R_n$ are in the
normal closure of the $R_i$ for $i < n$.

We consider how to construct $A=\lan \lan R_*\ran\ran$. In dimension
1, this is the normal closure of $R_1\cup \delta_2(R_2)$ as extended
to the groupoid case  in \cite{brownbook:2,Higginsbook}, i.e we take
the `consequences' of $R_1\cup \delta_2(R_2)$ in $C_1$. Suppose this
$A_1$ has been constructed.

\begin{proposition}\label{prop:normalclosure}
For $n \geq 2$, $A_n=\lan\lan R_* \ran\ran_n$ is generated as a
group (abelian if $n \geq 3$) by the elements
$$r^c, \; x-x^a \; \text{ for all }  r \in R_n \cup \delta_{n+1}(R_{n+1}),
\; c \in C_1, \; x \in C_n, \; a \in A_1.$$
\end{proposition}
\begin{proof}
Clearly these elements belong to $A_n$. We now prove the set
of these is $C_1$-invariant.

This is clear for the set of elements of the form $r^c$ as above.
Suppose then $x,c,a$ are as above. Then \begin{align*} (x-x^a)^c&=
x^c - x^{a+c}\\ &= x^c - (x^c)^{-c+a+c},
\end{align*}
which implies what we want since $-c+a+c \in A_1$ by normality.

It follows that the group generated by these elements is
$C_1$-invariant. In dimension 2, this implies normality, by the
crossed module rules.
\end{proof}

\section{Free crossed complexes}\label{sec:free}
Write $\F(n)$ for the crossed complex freely generated by one
generator $c_n$ in dimension $n$. So $\F(0)$ is a singleton in all
dimensions; $\F(1)$ is essentially the groupoid $\I$\,; and for $ n
\ge 2$,  $\F(n)$ is in dimensions $n$ and $n-1$ an infinite cyclic
group with generators $c_n$ and $\delta c_n$ respectively, and
otherwise trivial. Let $\bS(n-1)$ be the subcrossed complex of
$\F(n)$ generated by $\delta c_n$. Thus $\bS(-1)$ is empty.

If $E^n_*$ and $S^{n-1}_*$ denote the skeletal filtrations of the
standard $n$-ball and $(n-1)$-sphere respectively, then a basic
result in algebraic topology is that $$\Pi E^n _* \cong \F(n), \quad
\Pi S^{n-1}_* \cong \bS(n-1).$$ This is also a consequence of the
Higher Homotopy van Kampen Theorem  indicated in  \ref{thm:colim},
see \cite{bh:colimits}.

We now define a particular kind of morphism $j: A \to F$ of crossed
complexes which we call a {\it morphism of relative free type}. Let
$A$ be any crossed complex. A sequence of morphisms $j_n: F^{n-1}
\to F^{n}$ may be constructed inductively as follows. Set
$F^{-1}=A$. Supposing $F^{n-1}$ is given, choose any family of
morphisms as in the top row of the diagram
$$ \xybiglabels\xymatrix{\bigsqcup _{\lambda \in \Lambda_n}\mathbb{S}{(n-1)}
\ar[rr]^-{(f^{\lambda })} \ar[d] & & F^{n-1} \ar[d]^{j_{n}} \\
\bigsqcup _{\lambda \in \Lambda_n} \mathbb{F}(n) \ar[rr] & & F^n }$$
and form the pushout in $\Crs$ to obtain $j_n:F^{n-1}\to F^n$. Let
$F= \colim _n F^n$, and let $j: A \to F$ be the canonical morphism.
The image $x^n_\lambda$ of the element $c_n$ in the summand indexed
by $\lambda$ is called a {\it basis element of $F$ relative to $A$},
and we may conveniently write
$$F = A \cup \{ x^n_\lambda  \}_{\lambda \in \Lambda_n, n \geq 0}. $$

We now give some useful results on this notion, see
\cite{bh:class91}.
\begin{proposition}\label{propII1:freecomp}
Given two morphism of relative free type, so is their composite.
\end{proposition}

\begin{proposition} \label{propII1:freepo}If in a pushout square
$$\xybiglabels \xymatrix{A \ar[r] \ar[d] & A' \ar[d] \\
F \ar[r] & F'}$$ \noindent the morphism $A \to F$ is of relative
free type, so is the morphism $A' \to F'$.
\end{proposition}

\begin{proposition} \label{prop:freecolim}
If in a commutative diagram
$$\xybiglabels \xymatrix{A^0 \ar[r] \ar[d] & A^1 \ar[r] \ar[d] & \ldots \ar[r] & A^n \ar[r] \ar[d] & \ldots \\
F^0 \ar[r] & F^1 \ar[r] & \ldots \ar[r] & F^n \ar[r] & \ldots }$$
each vertical morphism is of relative free type, so is the induced
 morphism $\colim_n A^n \to \colim_n F^n$.
\end{proposition}

In particular:
\begin{corollary} \label{corII1:freeseq}
If in a sequence of morphisms of crossed complexes
 $$F^0 \to F^1 \to \cdots \to F^n \to \cdots$$
each morphism is of relative free type, so are the composites $F^0
\to  F^n$ and the induced morphism $F^0 \to \colim_n F^n$.
\end{corollary}

A crossed complex $F$ is {\it free on $R_*$} if in the first place
$R_*$ generates $F$, and secondly morphisms on $F$ to any crossed
complex can be defined inductively by their values on $R_*$.

So in the first instance we have $R_0= F_0$, and $F_1$ is the free
groupoid on the graph $(R_1,R_0,s,t)$. We assume this concept as
known; it is fully treated in \cite{brownbook:2,Higginsbook}

Secondly, $R_2$ comes with a function $w: R_2 \to F_1$ given by the
restriction of $\delta_2$. We require that the inclusion $R_2 \to
F_2$ makes $F_2$ the free crossed $F_1$-module on $R_2$.

By this stage, the fundamental groupoid $\pi_1 F$ is defined; we
require that for $n \geq 3$, $F_n$ is the free $\pi_1 F$-module on
$R_n$.

A standard fact, due in the reduced case to Whitehead in [CHII],  is
that if $X_*$ is the skeletal filtration of a $CW$-complex, then
$\Pi X_*$ is the free crossed complex on the characteristic maps of
the cell structure of $X_*$. This may be proved using the relative
Hurewicz theorem, and is also a consequence of the Higher Homotopy
van Kampen Theorem{\it Higher Homotopy  van Kampen Theorem,
(HHvKT)}\footnote{Jim Stasheff has recently suggested this term to
Brown should replace the previous  `Generalised van Kampen Theorem'
in order to emphasise the nature of the theorem. } of
\cite{bh:colimits}.

\begin{proposition} \label{prop:specifymorph}
If $C$ is a free crossed complex on {$R_*$}, then a morphism $ f: C
\to D$ is specified by the values $fx \in D_n,  x \in R_n, n \ge 0$
provided only that the following geometric conditions hold:
\begin{equation}sfx=fsx,
 x \in R_1, tfx=ftx, x \in R_n, n \ge 1, \delta fx
=f \delta x, x \in R_n, n \ge 2.  \end{equation}
\end{proposition}
We refer also to \cite{brown-gol,bh:class91} for more details on
free crossed complexes. It is proved in \cite{brown-gol} that a weak
equivalence of free crossed complexes is a homotopy equivalence.

We now illustrate some of the difficulties of working with free
crossed modules by giving a proposition and a counterexample due
essentially to Whitehead, \cite{wjhc:sht}.

\begin{theorem} \label{thm:free-sub} Let $C$ be the free crossed complex
on $R_*$,  and suppose $S_* \subseteq R_*$ generates a subcrossed
complex $B$ of $C$. Let $F$ be the free crossed complex on $S_*$.
Then the induced morphism $F \to C$ is injective if the induced
morphism $\pi_1 B \to \pi_1 C$ is injective.
\end{theorem}
\begin{proof}
First of all, we know that a subgroupoid of a free groupoid is free.
Also in dimensions $>2$ $C_n$ is the free $\pi_1 C$-module on the
basis $R_n$. So injectivity, under the given condition,  is clear in
this case.

Thus the only problem is in dimension 2, and here we generalise an
argument of Whitehead, \cite{wjhc:sht}. We use the functor $\nabla:
\Crs \to \Chn$ given in section
 \ref{sec:chains}.

The abelianised groupoids $F_2^{\rm ab}, C_2^{\rm ab}$  are
respectively free $\pi_1 F,\pi_1 C$-modules  on the bases $S_2,R_2$.
Since the induced morphism on $\pi_1$  is injective, so also is the
induced morphism $F_2^{\rm ab}\to C_2^{\rm ab}$. The morphism $C_2
\to C_2^{\rm ab}$ is injective on $\Ker \delta_2: C_2 \to C_1$,
since $C_1$ is a free groupoid,  by (ii) of
\ref{thmII1:crossedtochain}. So $F \to C$ is injective in dimension
2.
\end{proof}

\begin{example}
Let $X=Y= \{ x\}, R=\{a,b\}, S= \{b\}$ where $a=x, b=1$. The group
presentations  $\langle Y\mid S\rangle, \langle X\mid R \rangle$
determine free crossed modules $\delta_S:C(S) \to F(X),
\delta_R:C(R) \to F(X)$. The inclusion $i: S \to R$ determines
$C(i): C(S) \to C(R)$. Now $F(X)=F(Y)=\C$, the infinite cyclic
group, while $C(S)$ is abelian and is the free $\C$-module on the
generator $b$. Also in $C(R)$, $ab=ba$ since $\delta_R b=1$. Hence
$$C(i)(b^x)=(C(i)(b))^{\delta_R a}= a^{-1}ba = b=C(i)(b),$$ and so
$C(i)$ is not injective.

{Of course the geometry of this example is the cell complex $K= E^2
\vee S^2$ and the subcomplex $S^1 \vee S^2$.}  \hfill $\Box$
\end{example}

\section{Cylinder and homotopies} \label{sec:cyl-homotopy}
It is useful to write out first all the rules for the cylinder
$\mathrm{Cyl} \, (C) = {\mathcal I} \otimes C$, as a reference. For
full details of the tensor product, see \cite{bh:tens,brown:fields}.

Let $C$ be a crossed complex. The cylinder ${\mathcal I} \otimes
C$ is generated by elements $0 \otimes x, 1\otimes x$ of dimension
$n$ and $\iota \otimes x, (- \iota) \otimes x$ of dimension
$(n+1)$ for all $n \geqslant 0$ and $x \in C_n$,  with the
following defining
relations for $a\in \I$:\\
\noindent { \bf Source and target}
  \begin{alignat*}{2}
  t(a \otimes x) &=  t a \otimes t x &&
    \text{for all } ~ a \in {\mathcal I},  \in C\;\\
  s(a \otimes x) &= a \otimes s x
       && \text{if }~ a=0,1,n=1\;, \\
  s(a \otimes x) &= s a \otimes x
     && \text{if } a=\iota,- \iota,n=0\;.\\
\intertext{\bf Relations with operations} a \otimes x^{c} &= (a
\otimes x)^{t a \otimes c}\qquad  && \text{if }~ n \geqslant 2,\ c
\in C_1.\end{alignat*}
 \noindent {\bf Relations with additions}
\begin{align*} a \otimes (x+y) &=
\begin{cases}
(a \otimes x)^{t a \otimes y} + a \otimes y, & \text{if }
a=\iota,- \iota, n=1, \\
 a \otimes x + a \otimes y, &
\text{if } a=0,1, n \geqslant 1 \text{ or if } a=\iota, -\iota, \ n
\geqslant 2,
\end{cases} \\
(-\iota) \otimes x &= \begin{cases} -(\iota \otimes x)\hspace{6em} & \text{if } n=0,\\
-(\iota \otimes x)^{(-\iota) \otimes t x} &\text{if } n \geqslant
1.
\end{cases} \\
\intertext{\bf Boundaries } \delta(a \otimes x) & = \begin{cases}
a
\otimes \delta x & \text{if } a=0,1,\ n \geqslant 2;\\
-t a \otimes x - a \otimes s x + sa \otimes x + a \otimes
t x & \text{if } a=\iota, -\iota,\ n=1; \\
 -(a \otimes \delta x) - (ta \otimes x) + (sa \otimes x)^{a \otimes t x} &
\text{if } a=\iota,-\iota,\ n \geqslant 2.
\end{cases}
\end{align*}
\hfill $\Box$

Now we can translate the rules for a cylinder into rules for a
homotopy. Thus a {\it homotopy} $ f ^0 \simeq f $ of morphisms
$f^0, f : C \to D$ of crossed complexes is a pair $(h,f)$ where
$h$ is a family of functions $h_n : C_n \to D_{n+1} $ with the
following properties:
\begin{alignat}{2}
t h_n(x) &= t f(x) &\qquad &\mbox{for all } x \in C;\label{hom1} \\
h_1(x+y) &= h_1(x)^{f y}+  h_1(y)&\qquad &
\mbox{if } x,y \in C_1 \mbox{ and } x+y \mbox{ is defined;} \label{hom2}\\
 h_n(x+y) &= h_n(x) + h_n(y) &\qquad & \mbox{if } x,y \in
C_n,\; n \ge 2 \mbox{ and } x+y \mbox{ is defined;}\label{hom3} \\
h_n(x^{ c} )&= (h_n x)^{f c} &\qquad & \mbox{if }  x \in C_n,
 n \ge 2, \; c \in C_1, \mbox{ and } x^{c} \mbox{ is defined.} \label{hom4}
\end{alignat}  Then $f^0,f $ are related by
\begin{equation} f^0(x) =
\begin{cases}
s h_0 x & \mbox{if }  x \in C_0, \\
(h_0 s x)+(f x)+({\delta}_2h_1 x)-(h_0 t x) & \mbox{if }
 x \in C_1, \\
 \{ f x + h_{n-1}{\delta}_n x + {\delta}_{n+1}h_n x \}
^{-(h_0t x)} & \mbox{if }  x \in C_n, \; n \ge 2.
\end{cases}  \label{homreln} \end{equation}
\begin{remark} Part of the force of this statement is that if $(h,f)$
satisfy properties (11-14), then $f^0$ defined by (15) is a morphism
of crossed complexes.
\end{remark}

The following is a substantial result:
\begin{proposition}[\cite{bh:class91}] If $F,F'$ are free crossed
complexes, on bases $R_*,R'_*$,  then $F \otimes F'$ is the free
crossed complex on the basis $R \otimes R'$.
\end{proposition}
The proof in \cite{bh:class91} uses {the inductive construction of
free complexes as successive pushouts given in section
\ref{sec:free};  the exponential law and the symmetry of $\otimes$
show that $\otimes$ preserves colimits on either side, and this
gives an inductive proof, analogous to a corresponding result for
$CW$-complexes.}

A consequence, which may also be proved directly,  is:
\begin{proposition} \label{spechomot}
If $f:C \to D$ is a morphism of crossed complexes and $C$ is a free
crossed complex {on $R_*$}, then a homotopy $(h,f):f^0 \simeq f: C
\to D$ is specified by the values $hx \in D_{n+1}, x \in R_n, n \ge
0$ provided only that the following geometric conditions hold:
\begin{equation}
 t hx= t fx, x \in R_n, n \ge 0.
 \end{equation}
\end{proposition}
\begin{proof}
The main special fact we need here is that an $f$-derivation on a
free groupoid is uniquely defined by its values on a free basis.
This follows easily from the fact that an $f$-derivation  $h_1:C_1
\to D_2$ corresponds exactly to a section of a semidirect product
construction $F_1 \ltimes C_2 \to F_1$.
\end{proof}

\section{Cones and the HAL} \label{sec:cones-HAL}

\begin{Definition} \label{defn:cone}
Let $C$ be a crossed complex. The {\it cone } $\Cone(C)$ is defined
by the pushout
$$\xymatrix{\{1\} \otimes C \ar[d] \ar[r] & \{v\} \ar[d] \\
{\mathcal I} \otimes C \ar[r] & \Cone \,(C). }$$ We call $v$ the
{\it vertex} of the cone. \hfill $\Box$ \end{Definition}

Because the cone is formed from the cylinder by shrinking the end at
$1$ to a point, the rules for the  cylinder now simplify nicely.
\begin{proposition}\label{prop:Cone}
If $C$ is a crossed complex, then the cone  $\mathrm{Cone} \,(C)$
on $C$ is   generated by elements $0 \otimes x, \iota \otimes x, x
\in C_n$ of dimensions $n,n+1$ respectively, and $v$ of dimension
$0$
with the following rules, for all $a \in \I$:\\
\noindent { \bf Source and target} \begin{align*} t(a \otimes x)
&=
\begin{cases}
0 \otimes t x, & \text{if } a=0, \\
v & \text{otherwise. }\end{cases} \intertext{\bf Relations with
operations} a \otimes x^{c} &= a \otimes x \quad \text{if } \ n
\geqslant 2,\ c \in C_1.\\
\intertext{\bf Relations with additions} a \otimes (x+y) &= a
\otimes x + a \otimes y.\\ \intertext{and} (-\iota) \otimes x
&=\begin{cases}
-(\iota \otimes x) & \text{if } n=0,\\
-(\iota \otimes x)^{-\iota \otimes t x} & \text{if } n \geqslant
1.
\end{cases}\\
\intertext{\bf Boundaries } \delta_{n} (0 \otimes x) &= 0 \otimes
\delta_n x \qquad \text{if } n \geqslant 2.\\
 \delta_{n+1}(\iota \otimes x) & =
\begin{cases} - \iota \otimes s x + 0 \otimes x + \iota \otimes
t x & \text{if } n=1,\\  -(\iota \otimes \delta_n x) + (0 \otimes
x)^{\iota \otimes t x} & \text{if } n \geqslant 2.
 \end{cases}
\end{align*}
\end{proposition}
\begin{proposition}\label{prop:freecone}
Let $F$ be a free crossed complex on a basis $R_*$. Then
$\mathrm{Cone}(F)$ is the free crossed complex on $v$ in dimension
$0$, and elements $0 \otimes r, \iota \otimes r$ for all $r \in
R_*$, with boundaries given by proposition \ref{prop:Cone}.
\end{proposition}
\begin{proof}
This follows from proposition \ref{propII1:freepo}.
\end{proof}

We use the above to work out the fundamental crossed complex of the
simplex $a\Delta^n$ in an algebraic fashion. We define $a\Delta^0 =
\{0\}, \; a\Delta^n$ inductively by
$$a\Delta^n= \Cone(a\Delta^{n-1}). $$
The vertices of $a\Delta^1= \I$ are ordered as $0 < 1$. Inductively,
we get vertices $v_0, \ldots, v_n$ of $a\Delta^n$ with $v_n=v$ being
the last introduced in the cone construction, the other vertices
$v_i$ being $ (0, v_i)$. The fact that our algebraic formula
corresponds to the topological one follows from facts stated earlier
on the tensor product and on the GvKT.

We now define inductively top dimensional generators of the crossed
complex $a\Delta^n$ by, in the cone complex:
$$\sigma^0= v,\, \sigma^1=\iota, \, \sigma^n=( \iota \otimes
\sigma^{n-1}),\; n \geqslant 2,$$ with $\sigma ^0$ being the vertex
of $a\Delta^0$.

Next we need conventions for the faces of $\sigma^n$. We define
inductively
$$\partial_i \sigma^n= \begin{cases}
\iota \otimes \partial _i \sigma^{n-1} & \text {if } i < n ,\\
 0 \otimes \sigma ^{n-1} &\text { if } i = n.
\end{cases} $$

\begin{theorem}[Homotopy Addition Lemma] \label{thm:hal-simplex}
The following formulae hold, where $u_n= \iota \otimes v_{n-1}$:
\begin{align}
\delta_2 \sigma^2&=-\partial _1 \sigma^2 + \partial_2 \sigma^2 +
\partial_0 \sigma^2, \\
\delta_3 \sigma^3&= (\partial_3\sigma^3) ^{u_3} - \partial_0
\sigma^3 -\partial_2 \sigma^3 +\partial_1\sigma^3,\\
\intertext{while for $n \geqslant 4$} \delta_n\sigma^n& =(\partial_n
\sigma^n)^{u_n} + \sum _{i=0}^{n-1} (-1)^{n-i}
\partial_i\sigma^n .
\end{align}
\end{theorem}
\begin{proof}
For the case $n=2$ we have
\begin{align*}\delta_2 \sigma^2&=\delta_2 (\iota \otimes
\iota) \\
&= -\iota \otimes 0 + 0 \otimes \iota + \iota \otimes 1 \\
 &= -\partial _1 \sigma^2 + \partial_2 \sigma^2 + \partial_0
 \sigma^2.
 \end{align*}
 For $n=3$ we have:
\begin{align*}
 \delta_3 \sigma^3&=\delta_3 (\iota \otimes \sigma^2) \\
&= (0 \otimes \sigma^2)^{\iota \otimes v_2}
-\iota \otimes \delta_2 \sigma^2 \\
&= (0 \otimes \sigma^2) ^{u_3} - \iota \otimes (-\partial_1
\sigma^2+\partial_2 \sigma^2+
\partial_0 \sigma^2) \\
&= (\partial_3\sigma^3) ^{u_3} - \partial_0 \sigma^3 -\partial_2
\sigma^3 +\partial_1\sigma^3.\\ \intertext{We leave the general case
to the reader, using the inductive formula} \delta_{n+1}
\sigma^{n+1}&= (0 \otimes \sigma^n)^{\iota \otimes v_n} -\iota
\otimes \delta_n \sigma^n. \end{align*}  The key points that make it
easy are the rules on operations and additions of Proposition
\ref{prop:Cone}.
\end{proof}

\begin{corollary}
The formula for the boundary of a simplex is as given by the HAL in
section \ref{sec:hal}.
\end{corollary}
\begin{proof}
We use the fact that for $n \ge 2$, the geometric $n$-simplex
$\Delta^n$ may be regarded as the cone
$\mathrm{Cone}(\Delta^{n-1})$. Our previous results thus give an
isomorphism \begin{equation*} \Pi \Delta^n_* \cong \Cone(\Pi
\Delta^{n-1}_*).
\end{equation*}
Since $\Delta^1_* = E^1_*$, the HAL now follows from theorem
\ref{thm:hal-simplex}.
\end{proof}

\section{The unnormalised fundamental crossed complex \newline of a simplicial set}
\label{sec:PiK} We now give full details of the definition of the
(unnormalised) fundamental crossed complex of a simplicial set,
which we referred to in section \ref{sec:hal}.
\begin{Definition}
We define $\Piu K$ the {\it (unnormalised) fundamental crossed
complex of the simplicial set $K$} to be  the free crossed complex
having the elements of $K_n$ as generators in dimension $n$ and
boundary maps given by the Homotopy Addition Lemma. In detail this
gives the crossed complex $\Piu K$ as follows:
\begin{enumerate}[1.]
\item The objects are the vertices of $K$: $(\Piu K)_0 = K_0$;
\item The groupoid $(\Piu K)_1$ is the free groupoid associated to the
directed graph $K_1$. So it has a free generator $x: \partial_1 x
\to \partial _0 x$ for each $x \in K_1$;
\item The crossed module $(\Piu K)_2 \to (\Piu K)_1$ is the free
$(\Piu K)_1$-crossed module generated by the map $\delta_2 : K_2 \to
(\Piu K)_1$ given by
$$\delta_2x =  -\partial _1 x + \partial_2 x +
\partial_0 x$$ for all $x \in K_2$. We set $\pi_1 K= \Cok
\delta_2$, the {\it fundamental groupoid} of $K$.
\item For all $n \geqslant 3$, $(\Piu K)_n$ is the free
$\pi_1 K$-module with generators $K_n$ and boundary given by
\begin{equation*} \delta_n x  = \begin{cases} (\partial_3x) ^{u_3 x                                                                                               } -
\partial_0
x -\partial_2 x +\partial_1x & \text {if } n=3,\\
(\partial_n x)^{u_n x} + \sum _{i=0}^{n-1} (-1)^{n-i}
\partial_ix &  \text{if } n \ge 4,
\end{cases}
\end{equation*}
where $u_n= \partial^{n-1}_0$.
\end{enumerate}
 This construction is natural, giving a fundamental crossed
complex functor of simplicial sets  \begin{equation}\Piu : \Simp \to
\Crs. \tag*{$\Box$}\end{equation}
\end{Definition}
\begin{remark}
There are two notions of realisation of a simplicial set $K$,
usually written $\|K\|$, and $|K|$. In the first the only
identifications are along faces, and in the second the degenerate
simplices are also factored out. Each realisation is a $CW$-complex
with its skeletal filtration, and the Higher Homotopy van Kampen
Theorem of \ref{thm:colim}, \cite{bh:colimits}, shows that there is
a canonical isomorphism $\Piu K \cong \Pi (\|K\|_*)$.
\end{remark}
\section{0-normalisation}\label{sec:0-normalisation}
We first contrast with the usual case of a simplicial abelian group
$A$, where the simplicial operators $\partial_i, \eps_i$ are
morphisms of abelian groups. The {\it associated chain complex} $(A,
\partial) $ is then  $A_n$ in dimension $n\geq 0$ with boundary
$$\partial = \sum_{i=0}^n (-1)^i
\partial_i.$$ Let $(DA)_n$ for $n \geq 0$   be the subgroup of $A_n$
generated by the degenerate elements. It is  an easy calculation
from the rules for simplicial operators  that $\partial(DA)_n
\subset (DA)_{n-1}$ and so $(DA,\partial)$ is a subchain complex of
$(A, \partial)$.

In the nonabelian case, the formulae cope well with the increased
technicalities.  For the rest of this section, $K$ is a simplicial
set.
\begin{proposition}\label{prop:degen-normal}
Let $E_*$ be the set of degenerate elements in $ K$,  together with
the elements of $E_0$. Then $E_*$ is a normal structure in $\Piu K$.
\end{proposition}
\begin{proof}
By the rules $\partial_i \eps_i =\partial_{i+1} \eps_i= 1$, and the
Homotopy Addition Lemma,  we get immediate cancellation in $\delta_n
\varepsilon _i y$ for $0 < i < n-1$ but not necessarily for
$i=0,n-1$, because of the operators, and the nonabelian structures
in dimensions 1,2. Thus terms of concern are:
\begin{align*} \delta_2 \varepsilon _0 y &= -y+ \varepsilon _0
\partial_1 y +y,
\\
\delta_3 \varepsilon _0 y &= (\varepsilon _0 \partial_2
y)^{\partial_0 y}+ ( -y - \varepsilon _0 \partial_1 y +y), \\
\delta_3 \varepsilon _2 y &= (y)^{\varepsilon _0\partial_0^2 y}
-\varepsilon
_1 \partial_0y -y +\varepsilon _1 \partial_0 y,\\
&=(y^{\varepsilon _0\partial_0^2 y}-y) +(y   -\varepsilon _1
\partial_0y -y) +\varepsilon _1 \partial_0 y,\\  \intertext{and for
$n\geq 4$} \delta_n \varepsilon _{n-1} y&= y^{\varepsilon
_0\partial_0^{n-1}y}-y +\text{terms involving $\varepsilon _{n-2}
$}.
\end{align*}
This proves the result in view of the definitions in section
\ref{sec:normal}.
\end{proof}

\begin{Definition} We define a normal subcrossed complex
$E_0 K$ of $ \Piu K $ to be $K_0$ in dimension $0$ and in higher
dimensions to be normally generated by  degenerate elements of the
form $\varepsilon _0 y$.    \hfill $\Box$
\end{Definition}
\begin{Definition}
  We define the $0$-normalised crossed complex of
  $K$ to be $$\Pi ^{0N} K= (\Piu K)/E_0 K. $$
\end{Definition}
Our first result is:
\begin{theorem}
The projection $p^0: \Piu K \to \Pi ^{0N} K$ has a section $q$ such
that $qp^0 \simeq 1$.
\end{theorem}
The proof will occupy the rest of this section.

We first need a lemma, which will be used later as well.
\begin{lemma} \label{lem:deriv}
Let $h_1: (\Piu K)_1 \to (\Piu K)_2$ be a derivation. Then for $x
\in K_2$ we have\end{lemma}\vspace{-3ex}
\begin{align*} h_1 \delta_2 x&= -(h_1 \partial_1
x)^{\delta_2 x} +(h_1 \partial_2 x)^{\partial_0 x}+ h_1 \partial_0
x. \\[-3ex]
 \intertext{{\bf Proof }}
 h_1 \delta_2 x &=
h_1(-\partial _1 x +
\partial_2 x +
\partial_0 x) \\
&= (h_1 (-\partial _1 x + \partial_2 x))^{\partial_0 x} +
h_1\partial_0 x  \\
&= ((h_1 (-\partial_1 x))^{ \partial_2 x})^{\partial_0 x} + (h_1
\partial_2 x)^{\partial_0 x}+ h_1
\partial_0 x \notag
\\
&=-(h_1 \partial_1 x)^{\delta_2 x} +(h_1
\partial_2
x)^{\partial_0 x}+ h_1 \partial_0 x. \tag*{$\Box$}
\end{align*}
\begin{lemma}\label{lem:homdim1} If $h:\psi \simeq 1: \Piu K \to \Piu K$
is given by $h_0=\eps_0$ in dimension $0$, and in dimension $1$ by
$h_1$ is $\eps_0$ or $\eps_1$ on the free basis given by $K_1$, then
$\psi$ is given in dimensions $0,1$
by $$\psi x= \begin{cases}   x & \text{if } \dim x=0,\\
x-\eps_0\partial_0 x & \text{if } \dim x =1,
\end{cases} $$
and hence $\psi \eps_0 y= 0_y$ for all $y \in K_0$.
\end{lemma}
\noindent {\bf Proof} The case $\dim x =0$ is clear. For the case
$\dim x=1$ and for  $h_1=\eps_0$ we have
\begin{align*}
  \psi x &= \eps_0 sx +x +\delta_2(-\eps_0 x)-\eps_0 tx\\
          &= {\eps_0 \partial _1 x +x -(-x + \eps_0\partial_1 x+x ) -\eps_0 tx}\\
          &= x - \eps_0 \partial_0 x. \\
\intertext{and for  $h_1=\eps_1$ we have}  \psi x &= 0_{sx} +x +\delta_2(-\eps_1 x)-0_{tx}\\
            &= x +(x -x - \eps_0 \partial_0 x)\\
            &= x- \eps_0 \partial_0 x. \tag*{$\Box$}
\end{align*}

Now we define simultaneously a morphism $\psi: \Piu K \to \Piu K$
and a homotopy $h: \psi \simeq 1$ such that $\psi(E_0 K)$ is
trivial.

\begin{proposition}[$0$-normalisation]\label{normdim0} Let  $K$ be a simplicial set.
Then a  homotopy $(h,1)$ {on $\Piu K$ }may be defined on generators
from $K$ by $h_n=(-1)^n \varepsilon_0$, yielding $h: \psi \simeq 1$
where $\psi$ is given on generators by
$$\psi(x)= \begin{cases} x & \text{if } \dim x =0,\\
x-\varepsilon_0 \partial_0 x & \text{if } \dim x =1,\\
(x-\varepsilon_0 \partial_0 x)^{- \varepsilon_0 tx}& \text{if }
\dim x
>1.
\end{cases}$$ This $\psi$ satisfies

\noindent 1.-  $\psi(\varepsilon_0 x)=0_{tx}$ for all $x \in K$.

\noindent 2.- The induced {morphism  $\bar{\psi}:\Pi ^{0N} K \to
\Piu K$} satisfies $p_0 \bar{\psi}=1$ and $\psi = \bar{\psi }p_0
\simeq 1$. Thus $p_0$ is a homotopy equivalence.
\end{proposition}
\begin{proof}
To verify the formula for $\psi$ requires working out a formula for
$\epb_0 \delta_{n} x- \delta_{n+1}\varepsilon_0 x$, where $\epb_0$
is the derivation or operator morphism defined by $\eps_0$  on
generators, and we also have to use the crossed module rules.

Thus for $x \in K_2$, we have by Lemma \ref{lem:deriv}:
\begin{align*} \epb_0 \delta_2 x
&=-(\varepsilon_0 \partial_1 x)^{\delta_2 x} +(\varepsilon_0
\partial_2
x)^{\partial_0 x}+ \varepsilon_0 \partial_0 x\\ \intertext{while}
\delta_3 \varepsilon_0 x &=(\partial_3 \varepsilon_0
x)^{\partial^2_0\varepsilon_0 x} -x-\varepsilon_0\partial_1
  x+x \notag\\
  &= (\varepsilon_0 \partial_2 x)^{\partial_0x} +(-\varepsilon_0 \partial_1 x)^{\delta_2
  x}\notag \\
  &=(-\varepsilon_0 \partial_1 x)^{\delta_2
  x} +(\varepsilon_0 \partial_2 x)^{\partial_0x}, \quad \mbox{by centrality of $\delta_3 \varepsilon_0 x$}
\\ \intertext{From this we get}
-\delta_3\varepsilon_0 x + \epb_0 \delta_2 x &=\varepsilon_0
\partial_0 x.\\ \intertext{More easily, we have for $n \geqslant 3$
and $x \in K_n$} \delta_{n+1} \varepsilon_0 x &= (\varepsilon_0
\partial_n x)^{\partial^{n-1}_0x} +\sum
_{i=2}^{n}(-1)^{n+1-i}\partial_i \varepsilon_0 x \\
\intertext{and} \epb_0 \delta_n x &=(\varepsilon_0
\partial_n x)^{\partial^{n-1}_0x}+ \sum_{i=0}^{n-1} (-1)^{n-i}
\varepsilon_0
\partial_i x \\ \intertext{ so that}
\epb_0 \delta_{n} x- \delta_{n+1}\varepsilon_0 x &= (-1)^n
\varepsilon_0 \partial_0 x.
\end{align*}
With these computations we get $h : \psi \simeq 1$ where $\psi$ is
the morphism given in the statement. Hence $\psi(\varepsilon_0^n
v)=0_v$ for all $n \geqslant 1$, and in fact $\psi \varepsilon_0 x=
0_{tx}$ for all $x \in K$. From this we easily deduce that $\psi(\Pi
^0 K)$ is the trivial subcomplex on $K_0$.  The morphism $\psi$ then
defines a morphism $\bar{\psi}:\Piu ^{0N} K \to \Piu K$ satisfying
$\bar{\psi} p_0 = 1$.

The homotopy $\epb_0$ gives also $p_0 \bar{\psi} \simeq 1$. Thus
$\bar{\psi}$ is a homotopy equivalence (actually a deformation
retract). \phantom{x}
\end{proof}

\begin{remark} \label{rem:degvert} Let $v$ be a vertex of the simplicial set $K$. Then in
$\Piu K$ we have \begin{align*} \delta_2(\eps^2_0 v) &= \eps_0 v\\
\intertext{ and so $\eps_0v$ acts trivially on $(\Piu K)_n$ for $ n
\geq 3$. Further for $n \geq 3$} \delta_n (\eps^n_0 v) &=
\begin{cases} 0 & \text{ if } n \text{ is odd}, \\
\eps^{n-1}_0 v &\text{ if } n \text{ is even}. \tag*{$\Box$}
\end{cases}
\end{align*}
\end{remark}

\begin{proposition} The crossed complex $\Pi^{0N} K$ is isomorphic by $\bar{\psi}$ to
the (free) subcrossed complex of $\Piu K$ on the elements of $K$ not
of the form $\eps_0 y$ for $y \in K_{n-1}, \; n \geqslant 1$.
\end{proposition}
\begin{proof}
This follows from theorem \ref{thm:free-sub}.
\end{proof}
\begin{remark}
An advantage of working in the 0-normalised complex is that certain
awkward exponents, which would vanish or not appear in the usual
abelian case, now disappear in the 0-normalised complex. For example
if $y \in K_1$  we have \begin{alignat*}{2}   \delta_2 \eps_1y&=
-\partial_1 \eps_1 y +\partial_2 \eps_1 y + \partial_0 \eps_1 y
\qquad  & \delta_2 \eps_0y&= -\partial_1 \eps_0 y +\partial_2 \eps_0
y +
\partial_0 \eps_ y
\\&=-y+y + \eps_0 \partial _0 y  &
&= -y + \eps_0 \partial _1 y +y \\
&= 0_{ty}\qquad  \mod \eps_0. && = 0_{sy }\qquad  \mod \eps_0.
\tag*{$\Box$}
\end{alignat*}
\end{remark}
\begin{remark} There is another way of proceeding,  by first reducing
in $\Piu K$ all degeneracies of the vertices. Let $K_0$ denote also
the simplicial set on the vertices of $K$, and also the discrete
crossed complex on the object set $K_0$. Then the inclusion $K_0 \to
\Piu K_0$ is a strong deformation retract,  as is easily seen from
Remark \ref{rem:degvert}, with retraction $r_0: \Piu K \to K_0$,
say. So we may form the pushout
$$\xybiglabels \xymatrix{\Piu K_0 \ar [d] \ar [r] ^{r_0} & K_0 \ar [d] \\
\Piu K \ar [r] _r & \Pi^0 K}$$ Then $r$ is also a strong deformation
retract, by the methods of the homotopy theory of crossed complexes,
\cite{brown-gol}. We can then apply the previous methods to $\Pi^0
K$ to factor out the $0$-degeneracies. We leave details and
comparisons to the reader.
\end{remark}

\section{Normalisation} \label{sec:normalisation} Now we are able to define, in analogy
with Mac Lane \cite[\S VIII.6]{maclane:hom}, some further homotopies
on $\Pi^{0N}(K)$ to obtain the normalisation theorem.  We  can model
more closely the classical case on this $0$-normalised crossed
complex. Note that if $x \in K_n$ we write also $x$ for the
corresponding elements of both $\Piu K$ and $\Pi^{0N} K$.

\begin{Definition} For any $k \geqslant 0$ we define a subcrossed complex $D_kK \subset \Pi ^{0N} K
$ as follows:
    \begin{itemize}
    \item $(D_kK)_0 = (\Pi ^{0N} K)_0 =K_0$.
    \item $(D_kK)_1$ is trivial, i.e. consists only of identities.
    \item $(D_kK)_n$ is normally generated by $\eps_i y$ for  $y\in K_{n-1}$, $i \leqslant k$
    and $i \leqslant n-1$. %(equivalently $i \leqslant \min\{k, n-1\}$).
    \end{itemize}
Also, we define the {\it degeneracy subcomplex} $DK =
\bigcup_kD_kK$, i.e. $(DK)_n = \bigcup_k(D_kK)_n$ for all $n \in
\mathbb{N}$. \hfill $\Box$
\end{Definition}

Now we define a sequence of homotopies from the identity to
morphisms of crossed complexes sending $D_k K$ into $D_{k-1} K$ and
leaving fixed the elements up to dimension $k-1$. Then, the
composition of these morphisms is well defined and kills all the
degeneracy subcomplex. Let us formalise this sketch.

\begin{Definition} For any $k \geqslant 0$ we define a homotopy $(\tau^k,1):
\Pi ^{0N} K \to \Pi ^{0N} K$ given on the free basis  $x \in K_n$ by
\begin{align*}
\tau^k x &= \begin{cases} 0_{tx} &\text{if } n <k, \\
 (-1)^{n+k} \eps_k x & \text{if } n \geqslant k.
            \end{cases} \tag*{$\Box$}
\end{align*}
\end{Definition}
Therefore, for any $k\geqslant 0$ the homotopy $\tau^k$ defines a
morphism of crossed complex, $\phi^k: \Pi ^{0N} K \to \Pi ^{0N} K$
such that $\tau^k : \phi^k \simeq 1$. Clearly $\phi^0 = \psi $. For
$n \geqslant 1$ this map is given when $x \in K_n$ by
\begin{align*}
\phi^k x &= \begin{cases} x & \text{if }  n < k, \\
          x + (-1)^{k+n-1} \epb_k\delta_nx +
(-1)^{k+n}\delta_{n+1}\eps_kx &\text{if } k \leqslant n.
            \end{cases}
\end{align*}
where $\epb_i$ is the extension of $\eps_i$ on the basis to a
derivation or operator morphism as appropriate.

\begin{proposition} $\phi^k: \Pi ^{0N} K \to \Pi ^{0N}
K$ satisfies \begin{enumerate}[\rm (i)]
\item $\phi ^kD_jK \subset D_{j}K$ when $j<k$,  and
\item $\phi ^kD_kK \subset D_{k-1}K.$
\end{enumerate}
\end{proposition}
\begin{proof} \\ (i)  By the definition of $\phi^k$ we  have to prove the inclusion
only in the case $k \leqslant n$. In this case the generators of
$(D_jK)_n$ are elements $\eps _ix$ for $i\leqslant \min \{j,n-1\}$,
so the definition of $\phi^k$ is
$$\phi^k\eps_ix = \eps_ix + (-1)^{k+n-1}\eps_k\delta_n\eps_ix
+ (-1)^{k+n}\delta _{n+1}\eps_k\eps_i x.$$
 Therefore, since $\eps_ix
\in D_{j}K$, which is  a subcrossed complex, we have that $\delta
_n\eps_ix \in D_{j}K$. So  $\delta _n\eps_ix$ can be written as a
combination of $\eps_py$ with $y \in K_{n-2}$, $p\leqslant \min \{j,
n-2\}$. Therefore, since we have
$$\eps_k\eps_p =
\eps_p\eps_{k-1} \qquad \mbox{ if } k>p $$ we have that
$\eps_k\delta _n\eps_ix \in D_{j}K$.

On the other hand, for the same reason we have $\delta
_{n+1}\eps_k\eps_i\in D_{j}K$. Therefore, $\phi^k\eps_ix \in
D_{j}K$.

\noindent (ii) Now let us prove $\phi^kD_kK \subset D_{k-1}K$. Since
$(D_kK)_1$ is trivial we have to prove this inclusion only for
generators of dimension $n\geqslant 2$.

We first deal with the case $n=2$. Suppose then $x \in K_2$. Then
\begin{align*}
  \delta_3\eps_1x &= (\partial_3 \eps_1 x)^{\partial_0 x}
  -\eps_0\partial_0 x -x +x \\
  &=(\partial_3 \eps_1 x)^{\partial_0 x}  \qquad \mod \eps_0.\\
  \epb_1\delta_2 x&= (-\eps_1 \partial_1 x)^{\delta_2 x} +(\eps_1
  \partial_2 x)^{\partial_0 x} +\eps_1 \partial_0 x\\\intertext{so
  that $\mod \eps_0$ and by centrality}
  \phi^1 x &= x +\epb_1\delta_2 x - \delta_3\eps_1 x\\
  &= x - (\eps_1 \partial_1 x)^{\delta_2 x} +\eps_1\partial_0 x.
\end{align*}
Now it is clear that, $\mod \eps_0$, $x=\eps_1y $ implies $\phi^1
x=0$.

Let $\eps_i y \in (D_k K)_n$, where $i \leq \min\{k,n-1\}$. If $i
<k$ then $\eps_i y \in D_{k-1} K$ and so $\phi^k \eps_i y \in
D_{k-1} K$ by (i). It only remains to prove $\phi^k\eps_ky \in
D_{k-1}K$ for $y \in K_{n-1}$.

We have already done the case of $n \leq 2$.  In general
\begin{align*}\phi^k\eps_ky &= \eps_ky +
(-1)^{k+n-1}\eps_k\delta_n\eps_ky + (-1)^{k+n}\delta
_{n+1}\eps_k\eps_ky, \\ \intertext{for $y \in K_{n-1}$ with $n>2,$
and, in this case, $(D_{k-1}K)_n$ is abelian. We can write,}
\eps_k\delta_n\eps_ky  &= \eps_k (\partial _n\eps_ky)^{\partial
^{n-1}_0\eps_ky} + \sum_{j=0}^{n-1}(-1)^{n-j}\eps_k\partial
_j\eps_ky\\\intertext{ and}\delta_{n+1}\eps_k\eps_ky &= (\partial
_{n+1}\eps_k\eps_ky)^{
\partial_0^{n}\eps_k\eps_ky} +
\sum^{n}_{j=0}(-1)^{n+1-j}\partial_j\eps_k\eps_ky.
\end{align*}
Therefore $\phi ^k(D_kK) \subset D_{k-1}K$ follows from
\begin{align*} \eps_k\partial_j\eps_ky
&=\begin{cases}
\eps_{k-1}\eps_{k-1}\partial_jy & \mbox{if }  j<k\\
 \eps_{k}y & \mbox{if } j=k, k+1 \\
 \eps_{k}\eps_{k}\partial_{j-1}y & \mbox{if }
j>k+1\end{cases} \\ \intertext{and on the other hand,}
\partial_j\eps_k\eps_ky
& =\begin{cases}
  \eps_{k-1}\eps_{k-1}\partial_jy & \mbox{if } j<k \\
 \eps_{k}y & \mbox{if } j=k, k+1, k+2 \\
 \eps_{k}\eps_{k}\partial_{j-1}y & \mbox{if } j>k+2.
\end{cases}
\end{align*}
\end{proof}

Now we define $\phi = \phi _0 \phi^1 \cdots \phi^k \cdots : \Pi
^{0N} K \to \Pi ^{0N} K$.

Notice that since $\phi^k x = x$ for $k > \dim x$, this composite
is finite in each dimension.

\begin{proposition} $\phi DK = 0$.
\end{proposition}
\begin{proof} We have $(DK)_0 = 0$ and for $n >0$, $(DK)_n$ is generated by
$\eps_iy$ where $y \in K_{n-1}$ and $i \leqslant n-1$. Therefore,

$$\phi \eps_iy = \phi^0\phi^1 \cdots \phi^n\eps_iy$$

If $i = n-1$ we have that $\eps_iy \in (D_nK)_n$. So,
$$\phi^n\eps_iy \in D_{n-1}K, \quad \phi^{n-1}\phi^n\eps_iy \in
D_{n-2}K, \quad \cdots ,  \quad \phi^0 \cdots \phi ^n\eps_iy \in
D_0K .$$ If $i< n-1$ we have that $\eps_iy \in (D_iK)_n$. Therefore,
since $\phi^jD_iK \subset D_iK$ for $i<j$ we have $\phi^{i+1} \cdots
\phi^n\eps_iy \in D_iK$. So, as above , $\phi^0 \cdots \phi
^n\eps_iy \in D_0K $.
\end{proof}
\begin{Definition}
We define the {\it normalised fundamental crossed complex of the
simplicial set} $K$ by
$$\Pi  K = \frac{\Pi^{0N} K}{DK}.$$
\end{Definition}

\begin{theorem}The quotient morphism $$p : \Pi ^{0N} K \to \Pi K$$
is a homotopy equivalence with a section $q$. Further, $\Pi K$ has
free generators given by the images of the non degenerate elements
of $K$.
\end{theorem}
This follows as for the $0$-normalised case in the previous section.
Putting the two results together gives:
\begin{theorem} \label{thm:normalisation}The quotient morphism $$p : \Piu  K \to \Pi K$$
is a homotopy equivalence with a section $q$. Further, $\Pi K$ has
free generators given by the images of the non degenerate elements
of $K$.
\end{theorem}

%This last result
%is one place where cubical methods are more efficient, since for
%cubical sets $K$ and $L$ we have a cellular isomorphism $|K
%\otimes L| \cong |K| \times |L|$.

The crossed complex $\Pi  K$ homotopy equivalent to $\Piu K$ can be
described as freely generated by the non degenerate simplices of
$K$,  with boundary maps given by the HAL, forgetting the degenerate
parts. In this sense, we have two alternative descriptions of $\Pi
K$, one as just given and another  in terms of coends as described
in section \ref{sec:hal}. The latter is used in classifying space
results in \cite{bh:class91}.


\begin{thebibliography}{10}
\addcontentsline{toc}{section}{References} \expandafter\ifx\csname
url\endcsname\relax
  \def\url#1{{\tt #1}}\fi
\expandafter\ifx\csname urlprefix\endcsname\relax\def\urlprefix{URL
}\fi

\bibitem{baues1} H.~J. Baues, 1989, {\em Algebraic Homotopy\/}, volume~15 of {\em
Cambridge  Studies in Advanced Mathematics\/}, Cambridge Univ.
Press.


\bibitem{baues-br} H.-J. Baues and R. Brown, {\em On the relative homotopy groups of the
product  filtration and a formula of Hopf},  J. Pure Appl. Algebra
89 (1993) 49-61.

\bibitem{Baues-Cond} H.-J. Baues and  D. Conduch\'e, {\em  On the
tensor algebra of a nonabelian group and applications},  $K$-Theory
5 (1991/92) 531--554.

\bibitem{bauesT1}
H.-J. Baues and A. Tonks, {\em On the twisted cobar construction\/},
Math.   Proc. Cambridge Philos. Soc., 121 (1997) 229--245.



\bibitem{Blakers} A.~L. Blakers, {\em Some relations between homotopy
and homology groups}, Ann. Math. 49 (1948) 428--461.

\bibitem{breen} L. Breen, {\em On the classification of 2-gerbes and
2-stacks}, Ast\'erisque, 225 (1994) 1--160.

\bibitem{B67}
R.~Brown, \emph{Groupoids and {v}an {K}ampen's theorem}, {Proc.
{L}ondon Math. Soc.} (3) {17} (1967) 385--401.

\bibitem{B:coproducts} R. Brown,   \emph{Coproducts of crossed
$P$-modules: applications to second homotopy  groups and to the
homology of groups}, { Topology} 23 (1984) 337-345.

\bibitem{brownbook:2} R. Brown,  {\em Topology and Groupoids\/}, (2006)
Booksurge LLC, S. Carolina; (revised and retitled version of
previous editions: McGraw Hill, Maidenhead, 1968;  Ellis Horwood
Ltd., Chichester, 1988).

\bibitem{brown:hha}
R. Brown, {\em Groupoids and crossed objects in algebraic
topology\/},   Homology, homotopy and applications,  1 (1999) 1--78.

\bibitem{brown:fields}
R. Brown, {\em Crossed complexes and homotopy groupoids as non
commutative tools for higher dimensional local-to-global problems},
Proceedings of the Fields Institute Workshop on Categorical
Structures for Descent and Galois Theory, Hopf Algebras and
Semiabelian Categories, September 23-28, Fields Institute
Communications 43 (2004) 101-130.

\bibitem{brown:gilbert}
R. Brown and N.~D. Gilbert, {\em Algebraic models of 3-types and
automorphism  structures for crossed modules}, { Proc. London Math.
Soc.} (3) 59 (1989)  51-73.

\bibitem{brown-gol}R. Brown and   M. Golasinski, {\em A model
structure for the homotopy theory of crossed complexes}, { Cah. Top.
G\'eom. Diff. Cat}. 30 (1989) 61--82.

\bibitem{bh1978}
R. Brown and P.~J. Higgins, {\em On the connection between the
second relative  homotopy groups of some related spaces\/},
Proc.London Math. Soc., (3) 36 (1978) 193--212.


\bibitem{bh:algcub}
R. Brown and P.~J. Higgins, {\em The algebra of cubes\/}, J. Pure
Appl. Algebra,   21 (1981) 233--260.

\bibitem{bh:colimits}
R. Brown and P.~J. Higgins, {\em Colimit theorems for relative
homotopy   groups\/}, J. Pure Appl. Algebra, 22 (1981) 11--41.



\bibitem{bh:tens}
R. Brown and P.~J. Higgins, {\em Tensor products and homotopies for
  $\omega$-groupoids and crossed complexes,\/}, J. Pure Appl. Algebra, 4 (1987)
  1--33.

\bibitem{bh:chn}
R. Brown and P.~J. Higgins, {\em Crossed complexes and chain
complexes with   operators\/}, Math. Proc. Camb. Phil. Soc., 107
(1990) 33--57.

\bibitem{bh:class91}
R. Brown and P.~J. Higgins, {\em The classifying space of a crossed
complex\/}, Math. Proc. Camb. Phil. Soc. 110 (1991) 95-120.

\bibitem{Brown-Sivera} R. Brown, P.~J. Higgins  and R. Sivera, {\em Nonabelian
algebraic topology}, (in preparation). Part I downloadable.

\bibitem{brohu}R. Brown and J. Huebschmann, {\em Identities among relations}, in {\em
Low dimensional topology}, (ed. R.  Brown and T.L. Thickstun,
Cambridge University Press, 1982), London Math. Soc. Lecture Note
Series 48,   153-202.

\bibitem{brmopw}
R. Brown, E.~J. Moore, T. Porter and C.~D. Wensley, {\em Crossed
complexes, and free   crossed resolutions for amalgamated sums and
HNN-extensions of groups\/},   Georgian Math. J., 9 (2002) 623--644.


\bibitem{BRazak:LMS99}
R. Brown and A.~Razak Salleh, {\em Free crossed resolutions of
groups and presentations of modules of identities among
relations\/}, LMS J. Comput.  Math., 2 (1999) 28--61.

\bibitem{BS76}
R. Brown and C.~B. Spencer, {\em Double groupoids and crossed
modules\/},   Cahiers Topologie G\'eom. Diff\'erentielle, 17 (1976)
343--362.

\bibitem{BS1}
R. Brown and C.~B. Spencer, {\em $\mathcal G$-groupoids, crossed
modules and the  fundamental groupoid of a topological group\/},
Proc. Kon. Ned. Akad. v. Wet,  79 (1976) 296 -- 302.

\bibitem{br:w} R. Brown and C.~D. Wensley, {\em Computation and homotopical applications
of induced crossed modules}, J. Symbolic Computation 35 (2003)
59-72.

\bibitem{bullejosetal} M. Bullejos, E. Faro, and M.~A. Garcia-Mu\~{n}oz,
{\em Postnikov invariuants for crossed complexes}, J. Algebra, 285
(2005) 238-291.

\bibitem{crowell}
R. Crowell, {\em The derived module of a homomorphism\/}, Advances
in Math., 5   (1971) 210--238.

\bibitem{eil-mac-homhot} S. Eilenberg  and S. Mac Lane, {\em Relations between homology
and homotopy groups of spaces.  II}, Ann. of Math., (2) {51} (1950)
{514--533}.

\bibitem{eil-mac-acyclic} S. Eilenberg and S. Mac Lane, {\em Acyclic
models}, Amer. J. Math. {75} (1953)  {189--199}.

\bibitem{eil-mac} S. Eilenberg  and S. Mac Lane, {\em On the groups $H(\Pi,n)$:I},
Ann. of Math., (2) {58} (1953) {55-106}.

\bibitem{Higginsbook} P.~J. Higgins, {\em Categories and Groupoids},
van Nostrand, (1971), reprinted as a Theory and Applications of
Categories Reprint, 2005.

\bibitem{Hu-hal} S.-T. Hu,  {\em The homotopy addition theorem}, Ann. of Math. (2),
{58} (1953) 108--122.

\bibitem{Hueb} J. Huebschmann,  {\em Crossed $n$-fold extensions of groups and
cohomology},  Comment. Math. Helv. 55 (1980) 302--313.

\bibitem{Kamps-Porter} K.~H. Kamps and T. Porter, {\em Abstract homotopy and
simple homotopy theory}, World Scientific (1997).

\bibitem{kock} A. Kock, {\em Combinatorics of curvature and the
Bianci identity}, Theory and Applications of Categories, 2 (1996)
62-89.


\bibitem{maclane:hom}
S.~Mac Lane, 1967, {\em Homology\/}, number 114 in Grundlehren,
Springer.

\bibitem{maclane&jhcw}
S.~Mac Lane and J.~H.~C. Whitehead, {\em On the 3-type of a
complex\/}, Proc. Nat. Acad. Sci. U.S.A., 36, (1950), 41--48.

\bibitem{martins} J.~F. Martins, {\em On the Homotopy
Type and the Fundamental Crossed Complex of the Skeletal Filtration
of a CW-Complex}, Homology, Homotopy and Applications (to appear)
math.GT/0605364.


\bibitem{martins-porter} J.~F. Martins, T. Porter, {\em On Yetter's Invariant
and an Extension of the Dijkgraaf-Witten Invariant to Categorical
Groups},  math.QA/0608484.

\bibitem{emmathesis}
E.~J. Moore, 2001, {\em Graphs of Groups: Word Computations and Free
Crossed   Resolutions\/}, Ph.D. thesis, University of Wales, Bangor.


\bibitem{Peiff} R. Peiffer, {\em  \"Uber Identit\"aten zwischen Relationen},
Math. Ann., 121  (194)  67--99

\bibitem{Reid} K. Reidemeister, {\em  \"Uber Identit\"aten von Relationen},
Abh. Math. Sem. Universit\'at Hamburg, 16  (1949)  114--118.



\bibitem{sharko} V.~V. Sharko, {\em Functions on manifolds: algebraic and topological aspects},
Translations of Mathematical Monographs 131, American Mathematical
Society, Rhode Island (1993).




\bibitem{tonksthesis}
A.~P. Tonks, 1993, {\em Theory and applications of crossed
complexes\/}, Ph.D. thesis, University of Wales, Bangor.

\bibitem{GWWhitehead} G.~W. Whitehead, {\em Homotopy theory}, 1978
Springer.

\bibitem{jhcw:CHI}[CHI] J.~H.~C. {W}hitehead, {\em Combinatorial Homotopy {I}\/}, Bull.
Amer. Math.   Soc., 55 (1949) 213--245.

\bibitem{jhcw:CHII}[CHII]
J.~H.~C. {W}hitehead, {\em Combinatorial Homotopy II\/}, Bull. Amer.
Math.   Soc., 55 (1949) 453--496.

\bibitem{wjhc:sht} J.~H.~C. Whitehead, {\em Simple homotopy types\/}, Amer. J. Math.,
72 (1950)   1--57.

\end{thebibliography}
\end{document}